\documentclass[11pt]{article}
\usepackage{amsmath,amsfonts,amssymb,amsthm,bbm,mathptm,cases}
\usepackage{graphicx,graphics,epsfig,subfigure}
\usepackage{graphicx}
\usepackage{booktabs}
\usepackage{threeparttable}
\usepackage{xcolor,color,soul,array,cite,float,version,times}
\numberwithin{equation}{section}
\numberwithin{figure}{section}
\numberwithin{table}{section}
\numberwithin{footnote}{section}

\def\bx{\mathbf{x}}

\theoremstyle{definition}

\newtheorem{sch}{Scheme}[section]
\theoremstyle{plain}
\newtheorem{thm}{Theorem}[section]
\newtheorem{lem}{Lemma}[section]

\theoremstyle{remark}

\newcommand{\ben}{\begin{eqnarray}}
\newcommand{\een}{\end{eqnarray}}
\newcommand{\bea}{\begin{array}}
\newcommand{\eea}{\end{array}}
\newcommand{\bes}{\begin{subequations}}
\newcommand{\ees}{\end{subequations}}
\newcommand{\bef}{\begin{figure}[H]}
\newcommand{\eef}{\end{figure}}
\newcommand{\bet}{\begin{tikzpicture}}
\newcommand{\eet}{\end{tikzpicture}}

\def\bena#1\eena{\begin{eqnarray}\begin{array}{l}#1\end{array}\end{eqnarray}}
\def\besl#1\eesl{\begin{subequations}\begin{align}#1\end{align}\end{subequations}}

\newcommand{\parl}[2]{\ensuremath{\frac{\partial #1}{\partial #2}}}

\newcommand{\vparl}[2]{\ensuremath{\frac{\delta #1}{\delta #2}}}

\newcommand{\ohs}[1]{\ensuremath{\overline{#1}^{n+1/2}}}
\newcommand{\hs}[1]{\ensuremath{{#1}^{n+1/2}}}

\def\inc(#1){\includegraphics[height=3 cm]{pics/#1}}

\def\bx{\ensuremath{{\bf x}}}

\begin{document}
%\begin{CJK*}{GBK}{song}
%ÖÐÎIJâÊÔ
%\end{CJK*}
\title{ Error Estimates of Energy Stable Numerical Schemes for Allen-Cahn Equations with Nonlocal Constraints}
\author{
{Shouwen Sun} \footnote{sunshouwen@csrc.ac.cn, Beijing Computational Science Research Center, Beijing 100193, P. R. China. },
{Xiaobo Jing} \footnote{jingxb@csrc.ac.cn, Beijing Computational Science Research Center, Beijing 100193, P. R. China. }, and
%{Jia Zhao} \footnote{zhao62@math.sc.edu, Department of Mathematics, Utah State University, USA. }
{Qi Wang}\footnote{qwang@math.sc.edu, Department of Mathematics,
 University of South Carolina, Columbia, SC 29028, USA;  Beijing Computational Science Research Center, Beijing  100193, P. R. China; School of Materials Science and Engineering, Nankai University, Tianjin 300350, China .}
%%{Xiaogang Yang}\footnote{xgyang@csrc.ac.cn, Beijing Computational Science Research Center, Beijing  100193, P. R. China.}
%%{M. Gregory Forest}\footnote{forest@unc.edu, Departments of Mathematics and Biomedical Engineering, University of North Carolina at Chapel Hill,
%%Chapel Hill, NC 27599, USA.}
}
\date{\today}
\maketitle
\begin{abstract}
We present error estimates for four unconditionally  energy stable numerical schemes developed for solving Allen-Cahn equations with nonlocal constraints. The schemes are linear and second order in time and space, designed based on the energy quadratization (EQ) or the scalar auxiliary variable (SAV) method, respectively. In addition to the rigorous error estimates for each scheme, we also show that the linear systems resulting from the energy stable numerical schemes are all uniquely solvable.  Then, we present some numerical experiments to show the  accuracy of the schemes, their volume-preserving as well as energy dissipation properties in a drop merging simulation.

\end{abstract}
%%%%%%%%%%%%%%%%%%%%%%%%%%%%%%%%%%%%%%%%%%%%%%%%%%%%%%%%%%%%%%%%%%%%%%%%%%%%%%%%%%%%%%%%%%%%%%%%%%%%%%%%%%

\section{Introduction}

\noindent \indent The phase-field approach for multi-component material systems, in particular multi-phase fluids, has emerged as an important modeling and computational tool in recent years\cite{MORTON1996TWO, Du1991Numerical, Du2004A, Lowengrub1998Quasi, Teigen2011A, Yang2017Numerical, Yue2004A, Zhao2016Numerical,Boyer1999Mathematical}. The set of governing equations in phase field models for the material mixture is often derived variationally from its energy functional guided by the generalized Onsager principle or equivalently the second law of thermodynamics \cite{Yang2014Near}. This type of models is called thermodynamically consistent. Thermodynamically consistent models for material systems represent a class of partial differential equations that yield energy dissipation laws and possess variational structures. For thermodynamical models alone, which do not subject to the hydrodynamic constraints (mass, momentum, angular momentum, energy conservation),  these models are also known as gradient flow models, in which the time evolution of the thermodynamical variable is "proportional to" the variation of the system free energy. When the thermodynamic variable is a phase variable and the proportionality coefficient (mobility) is an algebraic operator, it's known as the  Allen-Cahn equation, whereas it is called a Cahn-Hilliard equation if the proportionality coefficient is a second order elliptic operator.

This class of models describes relaxation dynamics of the thermodynamical system to the equilibrium. In the case of a phase field description,  both the Allen-Cahn and Cahn-Hilliard models have been used to describe interfacial dynamics of multi-component immiscible fluids, where phase variables represent
volume fractions of material components. In these applications, the Cahn-Hilliard model conserves the phase volume for each material component while the Allen-Cahn does not warrant the conservation of the volume of the components. In the diffuse interface description of an otherwise fairly sharp interface, the use of the phase field model is primarily for maintaining a dynamical interface that allows topological evolution. In these cases, both the Cahn-Halliard and Allen-Cahn models are applicable. However, the lack of volume conservation in the Allen-Cahn model makes it less useful than the Cahn-Hilliard in the cases where volume is conserved. In order to conserve the volume in the Allen-Cahn model, the free energy functional has to be augmented by penalty terms that impose the volume conservation constraint or by a Lagrange multiplier to enforce the volume constraint explicitly \cite{Rubinstein1992Nonlocal, Li2017Unconditionally}. The Allen-Cahn models thus modified with a nonlocal penalizing constraint or a  nonlocal Lagrange multiplier constraint then become nonlocal equations.

Given their prominent role in describing diffuse interface dynamics and ease of use in many computations, the  Allen-Cahn and Cahn-Hilliard equation have been widely studied recently. Xu and Guo (cf. \cite{Guo2016High}) developed an LDG method for the Allen-Cahn equation and proved its energy stability.
Chen, Wang and Wise in \cite{Chen2014ALI} developed a linear iteration algorithm to implement a second-order energy stable
numerical scheme for a model of epitaxial thin film growth without slope selection.
Xu and Tang \cite{Xu2006Stability} used a different stabilizing mechanism to build stable large time-stepping, semi-implicit methods for an epitaxial growth model.
Lin presented a $C^0$ finite element method for a 2D hydrodynamic liquid crystal model which is simpler than existing $C^1$ element methods and mixed element formulation to preserve energy dissipation \cite{Lin2006Simulations}. A convex splitting strategy has been developed for phase field models in recent years by several research groups  \cite{Shen2012Second,Wang2011An,Wu2017Stabilized,Cai2017Error,Wang2012Unconditionally}. In the meantime, there is the linear stabilization approach \cite{Chen2015Decoupled,Shen2010Numerical,Zhao2016Energy,Yu2017Numerical} aimed for  developing energy stable numerical schemes by adding a high order stabilizing mechanism to the discrete scheme. More recently, the energy quadratization (EQ) method and its variant scalar auxiliary variable (SAV) method have simplified the development of energy stable schemes significantly and make their  development systematically \cite{yang2018linear,Yang2017Numerical2,zhao2017numerical}.

In the wake of the rapid development of various innovative strategies for energy stable approximations to dissipative partial differential equation systems,  there remains a serious task to obtain rigorous error analysis on the energy stable numerical schemes and their solvability.
In a recent work (cf. \cite{Guan2017Convergence}), Guan \emph{et al}. adopted the convex splitting approach to develop an energy stable scheme and obtained its error estimates, in which the convex part of the
nonlinear system is treated implicitly and the concave part explicitly in time marching schemes.
Yu and Wang in \cite{Wang2017Convergence} introduced the convergence analysis of an unconditionally energy stable linear Crank-Nicolson scheme for the Cahn-Hilliard Equation.
Yang in \cite{Yang2017Numerical1} considered  the Allen-Cahn and Cahn-Hilliard equation using the Invariant Energy Quadratization (IEQ)  method.
Wang and Wise  provide a detailed convergence analysis for an unconditionally
energy stable, second order accurate convex splitting scheme for the modified phase field crystal
equation \cite{Laird2012Convergence}.

In this paper, using the newly developed energy quadratization (EQ) technique and the scalar auxiliary variable (SAV) approach \cite{Yang2016Linear,Yang2017Numerical2,Yang2017Numerical,Zhao2017A}, we propose and analyze two pairs  of new second order schemes for the Allen-Cahn equation with nonlocal constraints, which consist of the Allen-Cahn model with a penalizing potential or a Lagrange multiplier that enforces the volume conservation for each material phase involved. In particular, we design the schemes so that they are  uniquely solvable and unconditionally energy stable. We present an extended class of schemes, numerical strategies, numerical implementation issues as well as performance tests of the models in another paper \cite{Jing&L&W2018}. In this one, we focus exclusively on error estimates for the four schemes. To the best of our knowledge, error estimates of the schemes for the Allen-Cahn equation with nonlocal volume conserving constraints  are not yet available. We present the error analyses on both the semi-discrete schemes as well as the fully discrete schemes. Since the proof for the fully discrete schemes are similar, we only detail the error estimate for a selected fully discrete scheme in this paper. Numerical tests for the accuracy of the schemes as well as one numerical example on merging drops are presented to illustrate the volume and energy conservation property as well as the usefulness of the schemes.

The rest of the paper consists of  the following sections. In \S 2, we present the derivation of the Allen-Cahn equation with a penalizing potential and a Lagrange multiplier, respectively. In \S 3,  we reformulate the models to equivalent forms via the EQ and SAV approach, respectively. We present the numerical schemes and state the energy dissipation property and the uniqueness of solution of the linear systems resulting from the schemes. In \S 4, we present rigorous error estimates for the schemes. In \S 5, we present  numerical convergence tests to demonstrate the accuracy of the schemes and test the volume-conservation and energy dissipation properties using a drop merging experiment. Finally, some concluding remarks are given in \S 6.

\section{Phase Field Models for Binary Material Systems}

\noindent \indent We consider a phase field model for time-dependent dynamics of binary material systems, in which phase A is represented by a phase variable at $\phi=0$ and phase B at $\phi=1$. We normally choose the phase variable $\phi\in[0,1]$ with $\phi$ identified as the volume fraction of material A and $1-\phi$ the volume fraction of material B in the binary material system.  The transition layer between the two phases is given by $0<\phi<1$ and "the interface" is defined heuristically at $\phi=\frac{1}{2}$ \cite{Gong2018Linear1}. The free energy of the system is denoted as $F(\phi)$, given by
\ben
F[\phi]=\int_\Omega[\frac{\gamma_1}{2}|\nabla\phi|^2+f(\phi)]{\mathrm{d\bx}},
\een
where $\gamma_1$ is a parameter measuring the strength of the conformational entropy and the bulk free energy density is denoted by $f(\phi)$, which is assumed bounded below, e.g., there exists a constant $C$ such that $f(\phi)\geq-C$.
In this paper, we will focus our study on a double well potential $f(\phi)=\gamma_2\phi^2(1-\phi)^2$, where $\gamma_2$ measures the strength of the double-well potential.

The simplest transport equation for time-dependent dynamics of $\phi$ is given by the Allen-Cahn equation
\ben
\bea{l}
\parl{\phi}{t}(\textbf{\emph{x}}, t)= -M \mu, \quad \textbf{\emph{x}}\in\Omega, t>0, \\
\parl{\phi}{n}(\textbf{\emph{x}} ,t)=0, \quad \textbf{\emph{x}}\in\partial\Omega, t>0,\\
\phi(\textbf{\emph{x}} ,t)|_{t=0}=\phi(\textbf{\emph{x}}, 0),\quad\textbf{\emph{x}}\in\Omega,
\eea\label{AC-eq}
\een
where $M$ is the mobility coefficient,  the chemical potential $\mu$ is given by
\bena
\mu=\vparl{F}{\phi}=\frac{\partial f}{\partial \phi}-\frac{\gamma_1}{2}\nabla\cdot\frac{\partial(|\nabla\phi|^2)}{\partial\nabla\phi}=-\gamma_1\Delta\phi+f'(\phi),
\eena
and $\Omega$ is the domain in which the binary material system occupies.
This model describes relaxation dynamics of the binary material system.

The energy dissipation rate of the Allen-Cahn equation is given by
\bena
\frac{d F}{dt}=\int_{\Omega} \vparl{F}{\phi}\phi_t {\mathrm{d\bx}}
=-\int_{\Omega} \mu (M \mu) {\mathrm{d\bx}} \leq 0,
\eena
provided   $M\geq 0$.

We denote the volume of material A by $V=\int_{\Omega} \phi {\mathrm{d\bx}}.$ Then,
\ben
\frac{dV}{dt}=-\int_{\Omega} M \mu {\mathrm{d\bx}}.
\een
It is normally not equal to zero, indicating that $V(t)$ is not conserved in the Allen-Cahn dynamics.

To apply this model to situations where volumes of the material components
are supposedly conserved, one has to modify it. One simple modification is to enforce the volume constraint $V(t)=V(0)$ by coupling it to the Allen-Cahn equation via a Lagrange multiplier or a penalizing potential.
The model so derived is called the   Allen-Cahn equation with nonlocal constraints.

\subsection {Allen-Cahn equation with nonlocal constraints}

\noindent \indent There exist a couple of ways one can impose the volume constraints to the Allen-Cahn dynamics. The first method is to penalize the difference between the volume and its initial value. We call the model the Allen-Cahn equation with a penalizing potential.

\subsubsection{Allen-Cahn model with a penalizing potential}

\noindent \indent The volume conservation means $V(t)=V(0)$ for any $t>0$. In this modified model, we   penalize $(V(t)-V(0))^2$   in the energy functional \cite{Li2017Unconditionally}.
Specifically, we modify  the free energy  as follows:
\bena
\tilde{F}=\int_{\Omega} [\frac{\gamma_1}{2}{|\nabla\phi|}^2+f(\phi)]{\mathrm{d\bx}}+\frac{\eta}{2}(\int_\Omega \phi(t){\mathrm{d\bx}}-V(0))^2,
\eena
where $\eta>0$ is a penalizing parameter, which is a model parameter and $V(0)=\int_\Omega \phi(\bx, 0) {\mathrm{d\bx}}.$

The transport equation for $\phi$ is given by \eqref{AC-eq}
with a modified chemical potential given by
\bena
\tilde \mu=\vparl{ \tilde{F}}{\phi}
=\mu+\sqrt{\eta}\zeta,
\eena
where
\bena
\zeta=\sqrt{\eta}(\int_\Omega \phi(t) {\mathrm{d\bx}}-V(0)).
\eena

The energy dissipation rate is given by
\bena
\frac{d F}{dt}=\int_{\Omega} \vparl{F}{\phi}\phi_t {\mathrm{d\bx}}
=-\int_{\Omega} \tilde \mu (M \tilde \mu) {\mathrm{d\bx}} \leq 0,
\eena
provided $M\geq 0$. The modified Allen-Cahn equation is approximately volume-preserving depending on the size of penalizing parameter $\eta>0$. In principle, the larger $\eta$ is, the more close $V(t)$ is to $V(0)$. However, in practice, this may not be the case. So, the choice of $\eta$ becomes an empirical trial.  We next discuss another approach to obtain a modified model that respects the volume conservation.

\subsubsection{Allen-Cahn model with a Lagrange multiplier}

\noindent \indent To enforce the volume conservation of each phase in the Allen-Cahn model, the model is modified  by strictly enforcing the volume constraint via a Lagrangian multiplier $L$ in the free energy \cite{Rubinstein1992Nonlocal}. Specifically,
 we modify the free energy by augmenting a penalty term with a Lagrange multiplier $L$ as follows
\bena
\tilde F=\int_{\Omega} [\frac{\gamma_1}{2}{|\nabla\phi|}^2+f(\phi)]{\mathrm{d\bx}}-L (V(t)-V(0)).
\eena

The transport equation in the modified Allen-Cahn equation for $\phi$ is given by \eqref{AC-eq}
with   the chemical potential given by
\bena
\tilde \mu=\vparl{ \tilde F}{\phi}
=\mu-L,
\eena
where
\bena
L=\frac{1}{\int_{\Omega}  M {\mathrm{d\bx}}}\int_{\Omega}[ M \mu]{\mathrm{d\bx}}.
\eena

This modified Allen-Cahn model now conserves the volume for each phase.
The energy dissipation rate is given by
\bena
\frac{d F}{dt}=\int_{\Omega} \vparl{\tilde F}{\phi}\phi_t{\mathrm{d\bx}}
=-\int_{\Omega} \tilde \mu (M \tilde \mu){\mathrm{d\bx}} \leq 0,
\eena
provided $M\geq 0.$.

The Allen-Cahn equations with nonlocal constraints not only conserve phase volume at least approximately in the case of a penalizing potential, but also dissipates energy. We next discuss how we approximate these equations numerically so that the properties are retained in the discrete solutions and the numerical solutions can be obtained efficiently.

\section{Energy Stable Numerical Approximations }

\noindent \indent For the  Allen-Cahn equations with nonlocal constraints, we employ the energy quadratization strategy and its variant scalar auxiliary variable approach developed recently to design numerical schemes to solve them, respectively.
The energy quadratization (EQ) and its variant, the scalar auxiliary variable (SAV) method, provide effective ways to derive linear, energy stable  numerical schemes \cite{Yang2016Linear,Yang2017Numerical2,Yang2017Numerical,Zhao2017A}. We discretize the equations in time using the linear Crank-Nicolson method firstly \cite{Wang2017Convergence}. For the nonlocal constraints, we use the Shermann-Morrison formula together with the composite Trapezoidal rule to ensure the second order accuracy and efficient implementation of the resulting linear system of equations.

Throughout the paper, we denote
\bena
\hs{f}=\frac{1}{2}(f^{n+1}+f^{n}),\quad
\overline{f}^{n+1/2}=\frac{3}{2}f^n-\frac{1}{2}f^{n-1}\\
\eena
 and let $\|\cdot\|$ denote the $L^2(\Omega)$ norm and $(,)$  the associated inner product of functions involved, respectively. All constants $C$ appearing in the paper represent generic constants independent of $\Delta t$.

\subsection{Numerical method for the Allen-Cahn model with a penalizing potential using EQ}

\noindent \indent In the   Allen-Cahn equation with a penalizing potential, we reformulate the free energy density by introducing two intermediate variables as follows
\bena
q=\sqrt{f(\phi)-\gamma_2 \phi^2+C_0},\quad
\zeta=\sqrt{\eta}(\int_\Omega \phi(t) \mathrm{d\bx}-V(0)).
\eena
Then, the free energy (2.6) recast to
\bena
\tilde{F}=\int_{\Omega} [\frac{\gamma_1}{2}|\nabla\phi|^2+\gamma_2{\phi}^2+q^2-C_0]\mathrm{d\bx}+\frac{\zeta^2}{2}.
\eena
The extended chemical potential is given by
 \bena
\mu=-\gamma_1 \Delta\phi+2 \gamma_2 \phi +2q g, \quad g=\parl{q}{\phi}, \\
\tilde u=\mu+\sqrt{\eta}\zeta.
 \eena
We reformulate the modified Allen-Cahn model as follows
\bena
\parl{\phi}{t}= -M \tilde \mu,\\
\parl{\zeta}{t}= \sqrt{\eta}\int_\Omega \parl\phi{t} \mathrm{d\bx},\\
\parl{q}{t}= g \phi_t.
\eena

Now we discretize it using the linear Crank-Nicolson method in time to arrive at a semi-discrete, second order numerical scheme.

\begin{sch} Given initial conditions $\phi^0$, we calculate $q^0$ from $\phi^0$ and  $\phi^1, q^1$ are calculated by a first order scheme. Having computed $\phi^{n-1}, q^{n-1}$ and $\phi^n, q^n$, we compute $\phi^{n+1},q^{n+1}$ as follows.
\bena
\phi^{n+1}-\phi^n= -\Delta t\ohs{M}  \hs{\tilde{\mu}},\\
\hs{\tilde{\mu}}=(\hs{-\gamma_1 \Delta\phi+2\gamma_2\phi)} +2\hs{q} \ohs{g}+\sqrt \eta \hs{\zeta},\\
\zeta^{n+1}-\zeta^{n}= \sqrt{\eta}\int_\Omega (\phi^{n+1}-\phi^n) \mathrm{d\bx},\\
q^{n+1}-q^n= \ohs{g}(\phi^{n+1}-\phi^n).
\eena
\end{sch}
We define the discrete energy as follows
\bena
F^n=\int_{\Omega} [\frac{\gamma_1}{2}(\nabla\phi^{n})^2+\gamma_2(\phi^{n})^2+({q}^{n})^2-C_0]\mathrm{d\bx}+\frac{(\zeta^n)^2}{2}.
\eena

This scheme is linear, second order in time, and the linear system resulting from it is uniquely solvable. The scheme obeys a discrete dissipation law, i.e., Scheme 3.1 is unconditionally energy stable\cite{Jing&L&W2018}, which implies $F^n\leq F^0$, namely,
$\|\phi^{n}\|_{L^2}, \|q^{n}\|_{L^2}, \|\zeta^{n}\|_{L^2} (n=1, 2, \ldots, N)$  are bounded.

\subsection{Numerical method for the Allen-Cahn model with a Lagrange multiplier using EQ}

\noindent \indent We reformulate the free energy density by introducing an intermediate variable
\bena
q=\sqrt{f(\phi)-\gamma_2{\phi}^2+C_0}.
\eena
Then, the free energy (2.10) recast to
\bena
\tilde{F}=\int_{\Omega} [\frac{\gamma_1}{2}|\nabla\phi|^2+\gamma_2{\phi}^2+q^2-C_0]\mathrm{d\bx}-L(V(t)-V(0)).
\eena
The extended chemical potentia is given by
\bena
\mu=-\gamma_1\Delta\phi+2 \gamma_2 \phi +2q g, \quad g=\parl{q}{\phi}, \\
\tilde u=u-L,\quad
L=\frac{1}{\int_{\Omega} M \mathrm{d\bx}}\int_{\Omega}[ M \mu]\mathrm{d\bx}.\label{Model3}
 \eena

We rewrite the modified Allen-Cahn equation with the nonlocal constraint as follows
\bena
\parl{\phi}{t}= -M \tilde \mu,\\
\parl{q}{t}= g \phi_t.
\eena
Using the linear Crank-Nicolson method in time, we obtain the following scheme.

\begin{sch} Given initial conditions $\phi^0$, we calculate $q^0$ from $\phi^0$ and   $\phi^1, q^1$ are computed using a first order scheme. Having computed $\phi^{n-1},q^{n-1}$, and $\phi^n,q^n$, we compute $\phi^{n+1},q^{n+1}$ as follows.
\bena
\phi^{n+1}-\phi^n= -\Delta t \overline{M}^{n+1/2} [\hs \mu-\hs{L}],\label {model12}\\
\hs \mu=(\hs{-\gamma_1 \Delta\phi+2\gamma_2\phi)} +2\hs{q} \ohs{g},\\
q^{n+1}-q^n= \ohs{g}(\phi^{n+1}-\phi^n),\label {model24}\\
\eena
where
\bena
\hs L=\frac{1}{\int_{\Omega} \ohs M \mathrm{d\bx}}\int_{\Omega}[ \ohs M \hs \mu]\mathrm{d\bx}.\\
\eena
\end{sch}
We define the discrete energy as follows
\bena
F^n=\int_{\Omega} [\frac{\gamma_1}{2}(\nabla\phi^{n})^2+\gamma_2{(\phi^{n})}^2+({q}^{n})^2-C_0]\mathrm{d\bx}.
\eena

Scheme 3.2 is linear, second order in time and unconditionally energy stable \cite{Jing&L&W2018}. The linear system resulting from the scheme is uniquely solvable. It implies  $F^n\leq F^0$, and
$\|\phi^{n}\|_{L^2}, \|q^{n}\|_{L^2} (n=1, 2, \ldots, N)$  are bounded.

The scalar auxiliary variable (SAV) method provides yet another strategy to arrive at linear, energy stable  numerical schemes.

\subsection{Numerical method for the  Allen-Cahn model with a penalizing potential using SAV}

\noindent \indent We now use the scalar auxiliary variable approach to design a new numerical scheme. We rewrite the energy function as follows
\bena
\tilde{F}=\int_{\Omega} [\frac{\gamma_1}{2}(\nabla\phi)^2+\gamma_2\phi^2]\mathrm{d\bx}+\int_{\Omega} [f(\phi)-\gamma_2\phi^2]\mathrm{\mathrm{d\bx}}+\frac{\eta}{2}(\int_\Omega \phi(t) \mathrm{\mathrm{d\bx}}-\int_{\Omega}\phi(0)\mathrm{\mathrm{d\bx}})^2.
\eena
We introduce an intermediate variable $E_1(\phi)=\int_{\Omega} [f(\phi)-\gamma_2\phi^2]\mathrm{\mathrm{d\bx}}$ and choose a constant $C_0$ such that $E_1(\phi)\geq -C_0$. We denote $U=\vparl{E_1}{\phi}$,  $r=\sqrt{E_1+C_0}$ (a scalar auxiliary variable) and $s(\phi)=\frac{\partial r}{\partial \phi}=\frac{U( \phi)}{2\sqrt{E_1( \phi)+C_0}}$, and define $\zeta=\sqrt{\eta}(\int_\Omega \phi(t) \mathrm{d\bx}-\int_{\Omega}\phi(0)\mathrm{\mathrm{d\bx}})$ as another scalar auxiliary variable. The free energy recast into
\bena
\tilde{F}=\int_{\Omega} [\frac{\gamma_1}{2}|\nabla\phi|^2+\gamma_2{\phi}^2]\mathrm{d\bx}+r^2-C_0+\frac{\zeta^2}{2}
\eena
and chemical potential is modified as follows
\bena
\tilde u=u+\sqrt{\eta}\zeta,
 \eena
where$$\mu=-\gamma_1 \Delta\phi+2 \gamma_2 \phi +2r s.$$

We reformulate the  modified Allen-Cahn equation with a penalizing potential as follows
\bena
\parl{\phi}{t}= -M \tilde \mu,\\
\parl{\zeta}{t}= \sqrt{\eta}\int_\Omega \parl\phi{t} \mathrm{d\bx},\\
\parl{r}{t}= \int_\Omega s\parl\phi{t} \mathrm{d\bx}.
\eena

Now we discretize it using the linear Crank-Nicolson method in time to arrive at a new semi-discrete scheme.
\begin{sch}
Given initial conditions $\phi^0$, we calculate $r^0$ and $\zeta^0$ from $\phi^0$ and then $\phi^1$,  $r^1$ and $\zeta^1$ are computed using a first order scheme. Having computed $\phi^{n-1}$, $r^{n-1}, \zeta^{n-1}$, $\phi^n$, $r^n$ and $\zeta^n$, we compute $\phi^{n+1}, r^{n+1}, \zeta^{n+1}$ as follows.
\bena
\phi^{n+1}-\phi^n= -\Delta t \overline{M}^{n+1/2} \hs{\tilde\mu},\\
\hs {\tilde\mu}=\hs{(-\gamma_1 \Delta\phi+2\gamma_2\phi)}+2\hs r \ohs s+\sqrt{\eta}\hs{\zeta},\\
\zeta^{n+1}-\zeta^{n}= \sqrt{\eta}\int_\Omega (\phi^{n+1}-\phi^n) {\mathrm{d\bx}},\\
r^{n+1}-r^n= \int_{\Omega} \ohs s(\phi^{n+1}-\phi^n){\mathrm{d\bx}},\\
\eena
where
\bena
\ohs s=\ohs {\frac{U( \phi)}{2\sqrt{E_1( \phi)+C_0}}}.\\
\eena
\end{sch}
We define the discrete energy as follows
\bena
F^n=\int_{\Omega}[\frac{\gamma_1}{2}(\nabla\phi^{n})^2+\gamma_2(\phi^{n})^2]\mathrm{\mathrm{d\bx}}+\frac{(\zeta^n)^2}{2}+({r}^{n})^2-C_0.
\eena

Scheme 3.3 is a linear, second order in time, unconditionally energy stable and the linear system resulting from the scheme is uniquely solvable \cite{Jing&L&W2018}. It then follows that the norms of solutions of the linear system resulting from Scheme 3.3:
$\|\phi^{n}\|_{L^2}$, $\|r^{n}\|_{L^2}$, $\|\zeta^{n}\|_{L^2} (n=1, 2, \ldots, N)$  are bounded.

\subsection{Numerical method for the Allen-Cahn  model with a Lagrange multiplier using SAV}

\noindent \indent We rewrite the energy functional as follows
\bena
\tilde{F}=\int_{\Omega} [\frac{\gamma_1}{2}(\nabla\phi)^2+\gamma_2\phi^2]\mathrm{d\bx}+\int_{\Omega} [f(\phi)-\gamma_2\phi^2]\mathrm{\mathrm{d\bx}}-L(\int_{\Omega}\phi(t)\mathrm{\mathrm{d\bx}}-\int_{\Omega}\phi(0)\mathrm{\mathrm{d\bx}}).
\eena
Introducing an intermediate variable  $E_1(\phi)=\int_{\Omega} [f(\phi)-\gamma_2\phi^2]\mathrm{\mathrm{d\bx}}$ and choosing a constant $C_0$ such that $E_1(\phi)\geq -C_0$, we denote $U=\vparl{E_1}{\phi}$ and $s(\phi)=\frac{\partial r}{\partial \phi}=\frac{U( \phi)}{2\sqrt{E_1( \phi)+C_0}}$. We then define $r=\sqrt{E_1+C_0}$ as the scalar auxiliary variable and recast the free energy into
\bena
\tilde{F}=\int_{\Omega} [\frac{\gamma_1}{2}|\nabla\phi|^2+\gamma_2{\phi}^2]\mathrm{d\bx}+r^2-C_0-L(V(t)-V(0))
\eena
and the extended chemical potential is given by
\bena
\tilde \mu=-\gamma_1\Delta\phi+2 \gamma_2 \phi +2r s-L, \\
L=\frac{1}{\int_{\Omega} M \mathrm{d\bx}}\int_{\Omega}[ M \mu]\mathrm{d\bx}.
\eena

We reformulate the  modified Allen-Cahn  model  as follows
\bena
\parl{\phi}{t}= -M \tilde \mu,\\
\parl{r}{t}= \int_\Omega s\phi_t \mathrm{d\bx}.
\eena

A new second order in time semi-discrete numerical scheme is obtained after applying the linear Crank-Nicolson method.
\begin{sch} Given initial conditions $\phi^0$, we calculate $r^0$ from $\phi^0$ and then  $\phi^1$ and $r^1$ are computed by a first order scheme. Having computed $\phi^{n-1}$, $r^{n-1}$, $\phi^n$ and $r^n$, we compute $\phi^{n+1}$ and $r^{n+1}$ as follows.

\bena
\phi^{n+1}-\phi^n= -\Delta t \overline{M}^{n+1/2} \hs{\tilde\mu},\\
\hs {\tilde\mu}=\hs {\mu}-\hs{L},\\
\hs {\mu}=(\hs{-\gamma_1 \Delta\phi+2\gamma_2\phi)} +2\hs{r} \ohs{s},\\
r^{n+1}-r^n= \int_{\Omega} \ohs s(\phi^{n+1}-\phi^n){\mathrm{d\bx}},\\
\eena
where
\bena
\ohs s=\ohs {\frac{U( \phi)}{2\sqrt{E_1( \phi)+C_0}}},\\
\hs L=\frac{1}{\int_{\Omega} \ohs M \mathrm{d\bx}}\int_{\Omega}[ \ohs M \hs \mu]\mathrm{d\bx}.
\eena
\end{sch}
We define the discrete energy as follows
\bena
F^n=\int_{\Omega}[\frac{\gamma_1}{2}(\nabla\phi^{n})^2+\gamma_2(\phi^{n})^2]\mathrm{\mathrm{d\bx}}+({r}^{n})^2-C_0,
\eena

Scheme 3.4 is a linear, second order in time, and unconditionally energy stable \cite{Jing&L&W2018}. The linear system resulting from the scheme is uniquely solvable.
 It follows that the norms of solutions of the linear system resulting from Scheme 3.4:
$\|\phi^{n}\|_{L^2}$, $\|r^{n}\|_{L^2}$ $(n=1, 2, \ldots, N)$ are bounded.

\subsection{{The fully discrete numerical approximations}}

\noindent \indent The semi-discrete numerical schemes have been proposed in the above subsections and the schemes are proven to be linear, second order in time, and unconditionally energy stable. In addition,  the linear systems resulting from the schemes are all uniquely solvable. Here, we consider the spatial discretization of the schemes to arrive at  fully discrete numerical schemes.
We use the finite difference method to discretize the  Allen-Cahn equations with nonlocal constraints in space.

Firstly, we begin with definitions of grid functions for the full discretization in two-dimensional space.
We denote $\Omega=[0,1]\times[0,1]$ as the computational domain and divide the domain into rectangular meshes with mesh size $h_x=1/N_x, h_y=1/N_y$, where $N_x$ and $N_y$ are two positive integers. We define the following 2D sets for grid points:
$$E_x=\{x_{i+\frac{1}{2}}|i=0,1,\ldots,N_x\}, C_x=\{x_{i}|i=1,2,\ldots,N_x\}, C_{\bar{x}}=\{x_{i}|i=0,1,\ldots,N_x+1\},$$
$$E_y=\{y_{j+\frac{1}{2}}|j=0,1,\ldots,N_y\}, C_y=\{y_{j}|j=1,2,\ldots,N_y\}, C_{\bar{y}}=\{y_{j}|j=0,1,\ldots,N_y+1\},$$
where $x_l=(l-\frac{1}{2})h_x, y_l=(l-\frac{1}{2})h_y$, $l$ can take on integer or half-integer values. In this paper, we chose $h_x=h_y=h$ for simplicity.

We define the following discrete function spaces
$$\mathcal{C}_{x\times y}=\{\phi: C_x \times C_y\rightarrow \mathbb{R}\}, \mathcal{C}_{\bar{x}\times y}=\{\phi: C_{\bar{x}} \times C_y\rightarrow \mathbb{R}\},$$
$$\mathcal{C}_{x\times \bar{y}}=\{\phi: C_x \times C_{\bar{y}}\rightarrow \mathbb{R}\}, \mathcal{C}_{{\bar{x}}\times \bar{y}}=\{\phi: C_{\bar{x}} \times C_{\bar{y}}\rightarrow \mathbb{R}\}.$$

The functions of $\mathcal{C}_{{\bar{x}}\times \bar{y}}$ are called cell centered functions. In component form, cell centered functions are identified via $\phi_{i,j}$. A discrete function $\phi\in \mathcal{C}_{{\bar{x}}\times \bar{y}}$ is said to satisfy homogeneous Neumann boundary conditions if and only if
\bena
\phi_{0,j}=\phi_{1,j}, \phi_{N_x,j}=\phi_{N_x+1,j}, j=1,2,\ldots,N_y. \phi_{i,0}=\phi_{i,1}, \phi_{i,N_y}=\phi_{i,N_y+1}, i=0,1,\ldots,N_x+1.
\eena

We define
$$D_x\phi_{i+\frac{1}{2},j}:=\frac{1}{h}(\phi_{i+1,j}-\phi_{i,j}), D_y\phi_{i,j+\frac{1}{2}}:=\frac{1}{h}(\phi_{i,j+1}-\phi_{i,j}),$$
and $\nabla_h, \Delta_h$ are denoted the discrete gradient operator, the discrete Laplace operator as follows
 $$\nabla_h \phi:=(D_x\phi, D_y\phi)^T, \Delta_h\phi:=\nabla_h \cdot \nabla_h \phi.$$

In addition, we define the following discrete inner products
$$(f,g):=h^2\sum_{i=1}^{M}\sum_{j=1}^N f_{i,j}g_{i,j},$$
$$[f,g]_x:=\frac{1}{2}h^2\sum_{i=1}^{M}\sum_{j=1}^N(f_{i+\frac{1}{2},j}g_{i+\frac{1}{2},j}+f_{i-\frac{1}{2},j}g_{i-\frac{1}{2},j}),$$
$$[f,g]_y:=\frac{1}{2}h^2\sum_{i=1}^{M}\sum_{j=1}^N(f_{i,j+\frac{1}{2}}g_{i,j+\frac{1}{2}}+f_{i,j-\frac{1}{2}}g_{i,j-\frac{1}{2}}).$$

For $\phi,\psi$, a natural discrete inner product of their gradients is given by
$$(\nabla_h \phi, \nabla_h \psi):=[D_x\phi, D_x\psi]_x+[D_y\phi, D_y\psi]_y,$$
and we also introduce discrete $\|\cdot\|_\infty$ norm, $\|\cdot\|_p$ norm ($1\leq p<\infty$) and discrete $\|\cdot\|_{H^1_h}$, $\|\cdot\|_{H^2_h}$,  $\|\cdot\|_{H^4_h}$ norm, respectively, as follows:
\bena
\|\phi\|_\infty:=max_{i,j}|\phi_{i,j}|,~~~
\|\phi\|_p:=(|\phi|^p,1)^{\frac{1}{p}},\\
\|\phi\|^2_{H^1_h}:=\|\phi\|^2_2+\|\nabla_h\phi\|^2_2,\\
\|\phi\|^2_{H^2_h}:=\|\phi\|^2_{H^1_h}+\|\Delta_h\phi\|^2_2,\\
\|\phi\|^2_{H^4_h}:=\|\phi\|^2_{H^2_h}+\|\nabla_h\Delta_h\phi\|^2_2+\|\Delta_h^2\phi\|^2_2.\\
\eena

\begin{lem}
For $\phi, \psi\in\mathcal{C}_{{\bar{x}}\times \bar{y}}$ satisfying the homogeneous Neumann boundary condition, the following summation by parts formulas can be derived:
\bena
-(\Delta_h \phi,\psi)=(\nabla_h \phi, \nabla_h \psi).
\eena
\end{lem}

 Secondly, we use the standard five-point method to discretize the Laplacian operator and apply the composite-trapezoid-formula to discretize integral terms to arrive at a second order scheme in space. On staggered grids in space, the composite-trapezoid-formula is given as follows
\bena
\int_\Omega f(\mathrm{{\bf x}})\mathrm{d{\bf x}}=h^2\sum_{i=1}^m\sum_{j=1}^nf({x_i, y_j})+O(h^2),
\eena
where $({x_i, y_j})\in C_x\times C_y.$
We denote the symbol
\bena
[1\star f]:=h^2\sum_{i=1}^m\sum_{j=1}^nf({x_i, y_j}).
\eena

Applying this spatial discretization strategy to the semi-discrete schemes developed earlier, we arrive at second order,  fully discrete numerical schemes.
Next, we only present the fully discrete scheme based on Scheme 3.2 as a representative as the others are discretized analogously and thus are omitted.

\begin{sch} Given initial conditions $\phi_{i,j}^0$, we calculate $q_{i,j}^0$ from $\phi_{i,j}^0$ and   $\phi_{i,j}^1, q_{i,j}^1$ are computed using a first order scheme. Having computed $\phi_{i,j}^{n-1},q_{i,j}^{n-1}$, and $\phi_{i,j}^n,q_{i,j}^n$, we compute $\phi_{i,j}^{n+1},q_{i,j}^{n+1}$ as follows.
\bena
\phi^{n+1}_{i,j}-\phi^n_{i,j}= -\Delta t \overline{M}^{n+1/2} [\hs \mu_{i,j}-\hs{L}],\\
\hs \mu_{i,j}=(-\gamma_1 \Delta_h\hs{\phi}_{i,j}+2\gamma_2\hs{\phi}_{i,j}) +2\hs{q}_{i,j} \ohs{g}_{i,j},\\
q^{n+1}_{i,j}-q^n_{i,j}= \ohs{g}_{i,j}(\phi^{n+1}_{i,j}-\phi^n_{i,j}),\\
\eena
where
\bena
\hs L=\frac{1}{\int_{\Omega} \ohs M \mathrm{d\bx}}[1\star (\ohs M \hs \mu)].\\
\eena
\end{sch}
The discrete energy is defined as follows
\bena
F^n=\frac{\gamma_1}{2}\|\nabla_h\phi^{n}\|^2_2+\gamma_2\|\phi^{n}\|_2^2+\|{q}^{n}\|^2_2-C_0|\Omega|.
\eena

Scheme 3.5 respects the dissipation law, and thus is unconditionally energy stable. It then implies  $F^n\leq F^0$, and
$\|\phi^{n}\|_{2}, \|q^{n}\|_{2} (n=1, 2, \ldots, N)$  are bounded.
 We next present error estimates for the four semi-discrete schemes and the one fully discrete scheme presented above.

\section{Error Estimates }

\noindent \indent We consider a finite time interval [0,T] and a domain $\Omega\subset R^2$, which is open, connected and bounded with sufficiently smooth boundary $\partial\Omega$ such that the following Sobolev inequalities hold
\bena
\|f\|^2_{L^4(\Omega)}\leq C(2,\Omega)\|f\|_{L^2(\Omega)}\|f\|_{H^1(\Omega)},\\
\|f\|^2_{L^\infty(\Omega)}\leq C(2,\Omega)\|f\|_{H^1(\Omega)}\|f\|_{H^2(\Omega)}.\label{Sobolev ineq 4.1}
\eena
We denote $\Delta t=\frac{T}{N}.$
Let $\phi^n, q^n, g^n, \zeta^n, r^n, s^n$ be the numerical solutions obtained from the above schemes, and $\phi(t_n)$, $q(t_n)$, $g(t_n)$, $\zeta(t_n)$, $r(t_n)$, $s(t_n)$ be the exact solutions of the Allen-Cahn equations with nonlocal constraints evaluated at the discrete time, we define the error functions for $n=0,1,2,\ldots,N$ as follows
\ben
e^n_\phi=\phi(t_n)-\phi^n,  e^n_q=q(t_n)-q^n, e^n_g=g(t_n)-g^n, e^n_\zeta=\zeta(t_n)-\zeta^n, e^n_r=r(t_n)-r^n, e^n_s=s(t_n)-s^n.
\een

For each scheme, we first establish the corresponding error equation and conduct a series of analyses to obtain the error estimate for the scheme. Without loss of generality, we set the mobility coefficient $M\equiv1$ in the following analysis.
\begin{lem}
Suppose there exists a positive constant $C$ such that
$$max_{n\leq N}(\|\phi(t_n)\|_{L^\infty}, \|\nabla\phi(t_n)\|_{L^\infty}, \|\phi^n\|_{L^\infty})\leq C,$$
then  the following inequalities hold
\bena
\|e^n_{g}\|_{L^2}\leq C\|e^n_\phi\|_{L^2},\\
\|\nabla e^n_{g}\|_{L^2}\leq C(\|e^n_\phi\|_{L^2}+\|\nabla e^n_\phi\|_{L^2}),\\
\|e^n_{s}\|_{L^2}\leq C\|e^n_\phi\|_{L^2},\\
\|\nabla e^n_{s}\|_{L^2}\leq C(\|e^n_\phi\|_{L^2}+\|\nabla e^n_\phi\|_{L^2}).\label{Lemma 4.1 inequality}
\eena
\end{lem}
\noindent {\bf Proof.} We know that the function $q(x)$ is $C^3$, $s(x)$ is $C^2$ for any $x\in\Omega$, then
\ben
|e^n_{g}|=|g(\phi(t_n))-g(\phi^n)|=|q^{''}(\xi)(\phi(t_n)-\phi^n)|\leq C |e^n_\phi|,
\een
which implies $\|e^n_{g}\|_{L^2}\leq C\|e^n_\phi\|_{L^2}$. Next, since $q^{'''}$ is continuous and $\nabla\phi(t_n)$ is bounded, we have
\begin{equation*}
\begin{aligned}
|\nabla e^n_{g}|&=|q^{''}(\phi(t_n))\nabla\phi(t_n)-q^{''}(\phi^n)\nabla\phi^n|\\
&=|(q^{''}(\phi(t_n))-q^{''}(\phi^n))\nabla\phi(t_n)+q^{''}(\phi^n)\nabla\phi(t_n)-q^{''}(\phi^n)\nabla\phi^n|\\
&=|q^{'''}(\xi)(\phi(t_n)-\phi^n)\nabla\phi(t_n)+q^{''}(\phi^n)\nabla e_\phi^n|\\
&\leq C(|e^n_\phi|+|\nabla e^n_\phi|).
\end{aligned}
\end{equation*}
where $\xi=\theta\phi(t_n)+(1+\theta)\phi^n$. The proof of other two inequalities are similar to the above and thus omitted.

\noindent
{\bf Assumption:} We assume that the exact solution $\phi$ of the Allen-Cahn equations with nonlocal constraints possesses the following regularity conditions
\bena
\phi\in L^\infty(0,T;H^2(\Omega))\bigcap L^\infty(0,T;W^{1,\infty}(\Omega)), \\
\phi_t\in L^\infty(0,T;L^\infty(\Omega))\bigcap L^2(0,T;H^1(\Omega)).\label{regularity estimation}
\eena

For the classical Allen-Cahn equation this is proved in \cite{karali2012existence}.

\subsection{Error estimates for the Allen-Cahn model of a penalizing potential  using EQ}

\noindent \indent In Scheme 3.1, the semi-discrete scheme is given as follows

\besl
&\frac{\phi^{n+1}-\phi^n}{\Delta{t}}-\hs{\gamma_1 \Delta \phi}+2 \gamma_2 \hs{\phi} +2\hs{q} \ohs{g}+\sqrt{\eta}\zeta^{n+1/2}=0,\label{Scheme3.1Numer Solu1}\\
&\frac{q^{n+1}-q^n}{\Delta t}= \ohs{g}\frac{(\phi^{n+1}-\phi^n)}{\Delta t},\label{Scheme3.1Numer Solu2}\\
&\frac{\zeta^{n+1}-\zeta^n}{\Delta t}=\sqrt{\eta}\int_\Omega\frac{\phi^{n+1}-\phi^n}{\Delta t}\mathrm{d\bx}.\label{Scheme3.1Numer Solu3}
\eesl

We substitute the exact solution of the modified Allen-Cahn equation and use the Taylor expansion to obtain

\besl
&\frac{\phi(t_{n+1})-\phi(t_{n})}{\Delta{t}}-\gamma_1\Delta\phi(t_{n+\frac{1}{2}})+2\gamma_2\phi(t_{n+\frac{1}{2}})
+2q(t_{n+\frac{1}{2}})\overline{g}(t_{n+\frac{1}{2}})+\sqrt{\eta}\zeta(t_{n+\frac{1}{2}})=R^{n+1/2}_\phi,\label{Scheme3.1Exact Solu1} \\
&\frac{q(t_{n+1})-q(t_{n})}{\Delta{t}}=\overline{g}(t_{n+\frac{1}{2}})\frac{\phi(t_{n+1})-\phi(t_{n})}{\Delta{t}}+R^{n+1/2}_q,\label{Scheme3.1Exact Solu2}\\
&\frac{\zeta(t_{n+1})-\zeta(t_{n})}{\Delta t}=\sqrt{\eta}\int_\Omega\frac{\phi(t_{n+1})-\phi(t_{n})}{\Delta{t}}\mathrm{d\bx}+R^{n+1/2}_\zeta,\label{Scheme3.1Exact Solu3}
\eesl
where $R^{n+1/2}_\phi$, $R^{n+1/2}_q$ ,$R^{n+1/2}_\zeta$ are corresponding truncation error and in the order of  ($O(\Delta{t}^2)$).

Subtracting (\ref{Scheme3.1Exact Solu1}), (\ref{Scheme3.1Exact Solu2}), (\ref{Scheme3.1Exact Solu3}) from (\ref{Scheme3.1Numer Solu1}), (\ref{Scheme3.1Numer Solu2}), (\ref{Scheme3.1Numer Solu3}), respectively, we arrive at  the error equations as follows

\besl
&\frac{e_\phi^{n+1}-e_\phi^{n}}{\Delta{t}}-\gamma_1\Delta{e}_\phi^{n+1/2}+2\gamma_2e_\phi^{n+1/2}
+2q^{n+1/2}{\overline{e}_{g}^{n+1/2}}+2e_q^{n+1/2}\overline{g}^{n+1/2}+\sqrt{\eta}e^{n+1/2}_\zeta=R^{n+1/2}_\phi,\label{error function 1} \\
&\frac{e_q^{n+1}-e_q^{n}}{\Delta{t}}={\overline{g}}^{n+1/2}\frac{e_\phi^{n+1}-e_\phi^{n}}{\Delta{t}}
+\overline{e}_{g}^{n+1/2}\frac{\phi(t_{n+1})-\phi(t_{n})}{\Delta{t}}+R_q^{n+1/2},\label{error function 2}\\
&\frac{e_\zeta^{n+1}-e_\zeta^{n}}{\Delta{t}}=\sqrt{\eta}\int_\Omega\frac{e_\phi^{n+1}-e_\phi^{n}}{\Delta{t}}\mathrm{d\bx}+R_\zeta^{n+1/2}.\label{error function 3}
\eesl

\begin{lem}
There exists a constant $C>0$ such that solution $\phi^n$ of Scheme 3.1 is bounded:
$$\|\phi^n\|_{L^\infty}\leq C,n=0,1,\ldots,N.$$
\end{lem}
\noindent {\bf Proof.} Firstly, we note that  $\|\phi^0\|_{L^\infty}\leq C$ is true by definition. We assume $\|\phi^n\|_{L^\infty} (n=1, 2, \ldots, k)(k< N)$ are bounded. Next, we use mathematical induction to prove $\|\phi^{k+1}\|_{L^\infty}$ is bounded.

We note that $q(\phi)$ and $g(\phi)$ are continuous as functions of  $\phi$ and $\phi^n(n=1, 2, \ldots, k)$ and are thus bounded.
Taking inner product of (\ref{error function 1}) with $e_\phi^{n+1}-e_\phi^n$ and $2\Delta{t}e^{n+1/2}_\phi$, respectively, we obtain
\bena
\frac{1}{\Delta{t}}\|e_\phi^{n+1}-e_\phi^{n}\|^2_{L^2}+\frac{\gamma_1}{2}(\|\nabla{e}_\phi^{n+1}\|^2_{L^2}-\|\nabla{e}_\phi^{n}\|^2_{L^2})
+\gamma_2(\|{e}_\phi^{n+1}\|^2_{L^2}-\|{e}_\phi^{n}\|^2_{L^2})+(2q^{n+1/2}{\overline{e}_{g}^{n+1/2}},e_\phi^{n+1}-e_\phi^n)\\
+(2e_q^{n+1/2}\overline{g(\phi)}^{n+1/2},e_\phi^{n+1}-e_\phi^n)+(\sqrt{\eta}e^{n+1/2}_\zeta,e_\phi^{n+1}-e_\phi^n)
=(R^{n+1/2}_\phi,e_\phi^{n+1}-e_\phi^n)\label{lemma4.2 1}
\eena
and
\bena
\|{e}_\phi^{n+1}\|^2_{L^2}-\|{e}_\phi^{n}\|^2_{L^2}+2\Delta{t}\gamma_1\|\nabla{e}^{n+1/2}_\phi\|^2_{L^2}+4\Delta{t}\gamma_2\|{e}^{n+1/2}_\phi\|^2_{L^2}+(2q^{n+1/2}{\overline{e}_{g}^{n+1/2}},2\Delta{t}e^{n+1/2}_\phi)\\
+(2e_q^{n+1/2}\overline{g(\phi)}^{n+1/2},2\Delta{t}e^{n+1/2}_\phi)+(\sqrt{\eta}e^{n+1/2}_\zeta,2\Delta{t}e^{n+1/2}_\phi)=(R^{n+1/2}_\phi,2\Delta{t}e^{n+1/2}_\phi).\label{lemma4.2 2}
\eena

Taking inner product of (\ref{error function 2}) with $2e_q^{n+1}$, we have
\bena
\|e^{n+1}_q\|^2_{L^2}-\|e^{n}_q\|^2_{L^2}+\|e^{n+1}_q-e^{n}_q\|^2_{L^2}\\
=({\overline{g}}^{n+1/2}(e_\phi^{n+1}-e_\phi^{n}),2e^{n+1}_q)+(\overline{e}_{g}^{n+1/2}(\phi(t_{n+1})-\phi(t_{n})),2e^{n+1}_q)
+\Delta{t}(R^{n+1/2}_q,2e^{n+1}_q).\label{lemma4.2 3}
\eena

Taking inner product of (\ref{error function 3}) with $2e_\zeta^{n+1}$, we obtain
\bena
\|e^{n+1}_\zeta\|^2_{L^2}-\|e^{n}_\zeta\|^2_{L^2}+\|e^{n+1}_\zeta-e^{n}_\zeta\|^2_{L^2}
=(\sqrt{\eta}\int_\Omega(e_\phi^{n+1}-e_\phi^{n})\mathrm{d\bx},2e^{n+1}_\zeta)+\Delta{t}(R^{n+1/2}_\zeta,2e^{n+1}_\zeta).\label{lemma4.2 4}
\eena

Combining (\ref{lemma4.2 1})-(\ref{lemma4.2 4}),we have
\bena
\|{e}_\phi^{n+1}\|^2_{L^2}-\|{e}_\phi^{n}\|^2_{L^2}+\gamma_2(\|{e}_\phi^{n+1}\|^2_{L^2}-\|{e}_\phi^{n}\|^2_{L^2})
+\frac{\gamma_1}{2}(\|\nabla{e}_\phi^{n+1}\|^2_{L^2}-\|\nabla{e}_\phi^{n}\|^2_{L^2})+\|e^{n+1}_q\|^2_{L^2}-\|e^{n}_q\|^2_{L^2}+\|e^{n+1}_\zeta\|^2_{L^2}\\-\|e^{n}_\zeta\|^2_{L^2}
+\frac{1}{\Delta{t}}\|e_\phi^{n+1}-e_\phi^{n}\|^2_{L^2}+\|e^{n+1}_q-e^{n}_q\|^2_{L^2}+\|e^{n+1}_\zeta-e^{n}_\zeta\|^2_{L^2}
+2\Delta{t}\gamma_1\|\nabla{e}^{n+1/2}_\phi\|^2_{L^2}+4\Delta{t}\gamma_2\|{e}^{n+1/2}_\phi\|^2_{L^2}\\
=({\overline{g(\phi)}}^{n+1/2}(e_\phi^{n+1}-e_\phi^{n}),2e^{n+1}_q)+(\overline{e}_{g}^{n+1/2}(\phi(t_{n+1})-\phi(t_{n})),2e^{n+1}_q)+(\sqrt{\eta}\int_\Omega(e_\phi^{n+1}-e_\phi^{n})\mathrm{d\bx},2e^{n+1}_\zeta)\\
-(2q^{n+1/2}{\overline{e}_{g}^{n+1/2}},e_\phi^{n+1}-e_\phi^n)-(2e_q^{n+1/2}\overline{g(\phi)}^{n+1/2},e_\phi^{n+1}-e_\phi^n)-(\sqrt{\eta}e^{n+1/2}_\zeta,e_\phi^{n+1}-e_\phi^n)\\
-(2q^{n+1/2}{\overline{e}_{g}^{n+1/2}},2\Delta{t}e^{n+1/2}_\phi)-(2e_q^{n+1/2}\overline{{g}(\phi)}^{n+1/2},2\Delta{t}e^{n+1/2}_\phi)-(\sqrt{\eta}e^{n+1/2}_\zeta,2\Delta{t}e^{n+1/2}_\phi)\\
+(R^{n+1/2}_\phi,e_\phi^{n+1}-e_\phi^n)+\Delta{t}(R^{n+1/2}_q,2e^{n+1}_q)+(R^{n+1/2}_\phi,2\Delta{t}e^{n+1/2}_\phi)+\Delta{t}(R^{n+1/2}_\zeta,2e^{n+1}_\zeta).\label{lemma4.2 5}
\eena

Using the Cauchy inequality with $\epsilon$, we obtain
\bena
|({\overline{g(\phi)}}^{n+1/2}(e_\phi^{n+1}-e_\phi^{n}),2e^{n+1}_q)|\\
=|((3{g}^n-{g}^{n-1})(e_\phi^{n+1}-e_\phi^{n}),e^{n+1}_q)|\\
\leq (\|3{g}^n\|_{L^\infty}+\|{g}^{n-1}\|_{L^\infty})|((e_\phi^{n+1}-e_\phi^{n}),e^{n+1}_q)|\leq \frac{C}{4\Delta{t}\epsilon}\|e_\phi^{n+1}-e_\phi^{n}\|^2_{L^2}+C\Delta{t}\epsilon\|e^{n+1}_q\|^2_{L^2}.\label{lemma 4.2 6}
\eena

From the General H\"{o}lder inequality, Sobolev inequalities (\ref{Sobolev ineq 4.1}) and inequality (\ref{Lemma 4.1 inequality}), and the regularity assumption (\ref{regularity estimation}), we obtain
\bena
|(\overline{e}_{g}^{n+1/2}(\phi(t_{n+1})-\phi(t_{n})),2e^{n+1}_q)|\\
=|((3e_{g}^n-e_{g}^{n-1})(\phi(t_{n+1})-\phi(t_{n})),e^{n+1}_q)|\\
\leq\|3e_{g}^n-e_{g}^{n-1}\|_{L^4}\|\phi(t_{n+1})-\phi(t_{n})\|_{L^4}\|e^{n+1}_q\|_{L^2}
\leq \|3e_{g}^n-e_{g}^{n-1}\|_{L^4}\|\phi'(t_{n}+\theta\Delta{t})\Delta{t}\|_{L^4}\|e^{n+1}_q\|_{L^2}\\
\leq C\Delta{t}\|3e_{g}^n-e_{g}^{n-1}\|^2_{L^4}+C\Delta{t}\|e^{n+1}_q\|^2_{L^2}
\leq C\Delta{t}(\|3e_{g}^n\|^2_{L^4}+\|e_{g}^{n-1}\|^2_{L^4})+C\Delta{t}\|e^{n+1}_q\|^2_{L^2}\\
\leq C\Delta{t}(\|e^n_\phi\|^2_{L^2}+\|\nabla{e}^n_\phi\|^2_{L^2}+\|e^{n-1}_\phi\|^2_{L^2}+\|\nabla{e}^{n-1}_\phi\|^2_{L^2})+C\Delta{t}\|e^{n+1}_q\|^2_{L^2}.\label{lemma 4.2 7}
\eena

By the General H\"{o}lder inequality, the Cauchy inequality with $\epsilon$, (\ref{Lemma 4.1 inequality}) and the boundedness of $\|q^{n}\|_{L^2}$ and $\|q^{n+1}\|_{L^2}$, we have
\bena
|(2q^{n+1/2}{\overline{e}_{g}^{n+1/2}},e_\phi^{n+1}-e_\phi^n)|\\
=|((q^{n+1}+q^{n})\frac{1}{2}(3e^n_{g}-e^{n-1}_{g}),(e_\phi^{n+1}-e_\phi^n))|\\
\leq \|3e_{g}^n-e_{g}^{n-1}\|_{L^4}\|q^{n+1}+q^{n}\|_{L^4}\|e_\phi^{n+1}-e_\phi^n\|_{L^2}
\leq \|3e_{g}^n-e_{g}^{n-1}\|_{L^4}(\|q^{n+1}\|_{L^2}+\|q^{n}\|_{L^2})\|e_\phi^{n+1}-e_\phi^n\|_{L^2}\\
\leq C\Delta{t}\epsilon(\|3e_{g}^n\|^2_{L^4}+\|e_{g}^{n-1}\|^2_{L^4})+\frac{C}{4\Delta{t}\epsilon}\|e_\phi^{n+1}-e_\phi^n\|^2_{L^2}\\
\leq C\Delta{t}\epsilon(\|e^n_\phi\|^2_{L^2}+\|\nabla{e}^n_\phi\|^2_{L^2}+\|e^{n-1}_\phi\|^2_{L^2}
+\|\nabla{e}^{n-1}_\phi\|^2_{L^2})+\frac{C}{4\Delta{t}\epsilon}\|e_\phi^{n+1}-e_\phi^n\|^2_{L^2}.\label{lemma 4.2 8}
\eena

Similarly, using the H\"{o}lder inequality, the Cauchy inequality with $\epsilon$, we arrive at the following estimates
\bena
|(2e_q^{n+1/2}\overline{{g}(\phi)}^{n+1/2},e_\phi^{n+1}-e_\phi^n)|\\
=|((e_q^{n+1}+e_q^{n})\frac{1}{2}(3{g}(\phi^n)-{g}(\phi^{n-1}),(e_\phi^{n+1}-e_\phi^n))|\\
\leq(\|3{g}(\phi^n)\|_{L^\infty}+\|{g}(\phi^{n-1})\|_{L^\infty})|(e_q^{n+1}+e_q^{n},e_\phi^{n+1}-e_\phi^n)|
\leq C\Delta{t}\epsilon(\|e^{n+1}_q\|^2_{L^2}+\|{e}^n_q\|^2_{L^2})+\frac{1}{4\Delta{t}\epsilon}\|e_\phi^{n+1}-e_\phi^n\|^2_{L^2},\label{lemma 4.2 9}
\eena
\bena
|(2q^{n+1/2}{\overline{e}_{g}^{n+1/2}},2\Delta{t}e^{n+1/2}_\phi)|\\
\leq \Delta{t}\|3e_{g}^n-e_{g}^{n-1}\|_{L^4}\|q^{n+1}+q^{n}\|_{L^4}\|e^{n+1/2}_\phi\|_{L^2}
\leq C\Delta{t}(\|3e_{g}^n\|^2_{L^4}+\|e_{g}^{n-1}\|^2_{L^4})+C\Delta{t}\|e^{n+1/2}_\phi\|^2_{L^2}\\
\leq C\Delta{t}(\|e^n_\phi\|^2_{L^2}+\|\nabla{e}^n_\phi\|^2_{L^2}+\|e^{n-1}_\phi\|^2_{L^2}+
\|\nabla{e}^{n-1}_\phi\|^2_{L^2})+\Delta{t}(\|e_\phi^{n+1}\|^2_{L^2}+\|e_\phi^{n}\|^2_{L^2}),\label{lemma 4.2 10}
\eena
and
\bena
|(2e_q^{n+1/2}\overline{{g}(\phi)}^{n+1/2},2\Delta{t}e^{n+1/2}_\phi)|
=\Delta{t}|((e_q^{n+1}+e_q^n)(3{g}(\phi^n)-{g}(\phi^{n-1})),(e^n_\phi+e^{n+1}_\phi))|\\
\leq\Delta{t}(\|3{g}(\phi^n)\|_{L^\infty}+\|{g}(\phi^{n-1})\|_{L^\infty})|(e_q^{n+1}+e_q^{n},e^n_\phi+e^{n+1}_\phi)|\\
\leq C\Delta{t}(\|e^{n+1}_q\|^2_{L^2}+\|{e}^n_q\|^2_{L^2})+C\Delta{t}(\|e_\phi^{n}\|^2_{L^2}+\|e_\phi^{n+1}\|^2_{L^2}).\label{lemma 4.2 11}
\eena

By applying the same method, we obtain the following results

\begin{equation}
\begin{aligned}
&|(R^{n+1/2}_\phi,e_\phi^{n+1}-e_\phi^n)|\leq\frac{1}{4\Delta{t}\epsilon}\|e_\phi^{n+1}-e_\phi^n\|^2_{L^2}+\epsilon\Delta{t}\|R^{n+1/2}_\phi\|^2_{L^2},\\
&|\Delta{t}(R^{n+1/2}_q,2e^{n+1}_q)|\leq\Delta{t}\|R^{n+1/2}_q\|^2_{L^2}+\Delta{t}\|e^{n+1}_q\|^2_{L^2},\\
&|(R^{n+1/2}_\phi,2\Delta{t}e^{n+1/2}_\phi)|\leq\Delta{t}\|R^{n+1/2}_\phi\|^2_{L^2}+\Delta{t}(\|e^{n+1}_\phi\|^2_{L^2}+\|e^{n}_\phi\|^2_{L^2}),\\
&|(\sqrt{\eta}e^{n+1/2}_\zeta,e_\phi^{n+1}-e_\phi^n)|
\leq \frac{C}{4\Delta{t}\epsilon}\|e_\phi^{n+1}-e_\phi^{n}\|^2_{L^2}+\Delta{t}\epsilon(\|e^{n+1}_\zeta\|^2_{L^2}+\|e^{n}_\zeta\|^2_{L^2}),\\
&|(\sqrt{\eta}\int_\Omega{(e_\phi^{n+1}-e_\phi^{n})}\mathrm{d\bx},2e^{n+1}_\zeta)|\leq\frac{C}{4\Delta{t}\epsilon}\|e_\phi^{n+1}-e_\phi^n\|^2_{L^2}+C\Delta{t}\epsilon\|e_\zeta^{n+1}\|^2_{L^2},\\
&|(\sqrt{\eta}e^{n+\frac{1}{2}}_\zeta,2\Delta{t}e^{n+1/2}_\phi)|\leq\Delta{t}(\|e_\zeta^{n+1}\|^2_{L^2}+\|e_\zeta^{n}\|^2_{L^2})
+\Delta{t}(\|e^{n+1}_\phi\|^2_{L^2}+\|e^{n}_\phi\|^2_{L^2}),\\
&|\Delta{t}(R^{n+1/2}_\zeta,2e^{n+1}_\zeta)|\leq\Delta{t}\|R^{n+1/2}_\zeta\|^2_{L^2}+\Delta{t}\|e^{n+1}_\zeta\|^2_{L^2}.\label{lemma 4.2 18}
\end{aligned}
\end{equation}

Choosing appropriate $\epsilon$ and combining the above inequalities (\ref{lemma 4.2 6})-(\ref{lemma 4.2 18}) with (\ref{lemma4.2 5}), we have
\bena
\|e^{n+1}_\phi\|^2_{L^2}+\gamma_2\|e^{n+1}_\phi\|^2_{L^2}+\frac{\gamma_1}{2}\|\nabla{e}^{n+1}_\phi\|^2_{L^2}+\|e^{n+1}_q\|^2_{L^2}+\|e^{n+1}_\zeta\|^2_{L^2}\\
-(\|e^{n}_\phi\|^2_{L^2}+\gamma_2\|e^{n}_\phi\|^2_{L^2}+\frac{\gamma_1}{2}\|\nabla{e}^{n}_\phi\|^2_{L^2}+\|e^{n}_q\|^2_{L^2}+\|e^{n}_\zeta\|^2_{L^2})\\
\leq C\Delta{t}(\|e^{n+1}_\phi\|^2_{L^2}+\|e^n_\phi\|^2_{L^2}+\|e^{n-1}_\phi\|^2_{L^2}+\|\nabla{e}^{n}_\phi\|^2_{L^2}+\|\nabla{e}^{n-1}_\phi\|^2_{L^2}
+\|e^{n+1}_q\|^2_{L^2}+\|e^n_q\|^2_{L^2}\\
+\|e^{n+1}_\zeta\|^2_{L^2}+\|e^n_\zeta\|^2_{L^2})
+\Delta{t}(\|R^{n+1/2}_\phi\|^2_{L^2}+\|R^{n+1/2}_q\|^2_{L^2}+\|R^{n+1/2}_\zeta\|^2_{L^2})\\
\leq C\Delta{t}(\|e^{n+1}_\phi\|^2_{L^2}+\|e^n_\phi\|^2_{L^2}+\|e^{n-1}_\phi\|^2_{L^2}\\
+\|\nabla{e}^{n+1}_\phi\|^2_{L^2}+\|\nabla{e}^{n}_\phi\|^2_{L^2}+\|\nabla{e}^{n-1}_\phi\|^2_{L^2}
+\|e^{n+1}_q\|^2_{L^2}+\|e^n_q\|^2_{L^2}\\+\|e^{n+1}_\zeta\|^2_{L^2}+\|e^n_\zeta\|^2_{L^2})
+\Delta{t}(\|R^{n+1/2}_\phi\|^2_{L^2}+\|R^{n+1/2}_q\|^2_{L^2}+\|R^{n+1/2}_\zeta\|^2_{L^2}).\label{lemma 4.2 19}
\eena

Summing (\ref{lemma 4.2 19}) over for time steps 1 to m ,we deduce that for $1\leq m\leq k$,
\bena
\|e^{m+1}_\phi\|^2_{L^2}+\frac{\gamma_1}{2}\|\nabla{e}^{m+1}_\phi\|^2_{L^2}+\|e^{m+1}_q\|^2_{L^2}+\|e^{m+1}_\zeta\|^2_{L^2}\\
\leq\|e^{m+1}_\phi\|^2_{L^2}+\gamma_2\|e^{n+1}_\phi\|^2_{L^2}+\frac{\gamma_1}{2}\|\nabla{e}^{m+1}_\phi\|^2_{L^2}+\|e^{m+1}_q\|^2_{L^2}+\|e^{m+1}_\zeta\|^2_{L^2}\\
\leq C\Delta{t}\sum^m_{n=0}(\|e^{n+1}_\phi\|^2_{L^2}+\|\nabla{e}^{n+1}_\phi\|^2_{L^2}+\|e^{n+1}_q\|^2_{L^2}+\|e^{n+1}_\zeta\|^2_{L^2})+N\Delta t\Delta{t}^4,\label{lemma 4.2 20}
\eena
and
\bena
min(1,\frac{\gamma_1}{2})(\|e^{m}_\phi\|^2_{L^2}+\|\nabla{e}^{m}_\phi\|^2_{L^2}+\|e^{m}_q\|^2_{L^2}+\|e^{m}_\zeta\|^2_{L^2})\leq\|e^{m}_\phi\|^2_{L^2}
+\frac{\gamma_1}{2}\|\nabla{e}^{m}_\phi\|^2_{L^2}+\|e^{m}_q\|^2_{L^2}+\|e^{m}_\zeta\|^2_{L^2}.\label{lemma 4.2 21}
\eena
Then
\bena
\|e^{m+1}_\phi\|^2_{L^2}+\|\nabla{e}^{m+1}_\phi\|^2_{L^2}+\|e^{m+1}_q\|^2_{L^2}+\|e^{m+1}_\zeta\|^2_{L^2}\\
\leq C\Delta{t}\sum^m_{n=0}(\|e^{n+1}_\phi\|^2_{L^2}+\|\nabla{e}^{n+1}_\phi\|^2_{L^2}+\|e^{n+1}_q\|^2_{L^2}+\|e^{n+1}_\zeta\|^2_{L^2})+N\Delta t\Delta{t}^4.\label{lemma 4.2 22}
\eena

Combining (\ref{lemma 4.2 22}) with the Gronwall's inequality, we have
\bena
\|e^{m+1}_\phi\|^2_{L^2}+\|\nabla{e}^{m+1}_\phi\|^2_{L^2}+\|e^{m+1}_q\|^2_{L^2}+\|e^{m+1}_\zeta\|^2_{L^2}\leq C\Delta{t}^4,\label{lemma 4.2 23}
\eena
which implies
\bena
\|e^{k+1}_\phi\|^2_{H^1}+\|e^{k+1}_q\|^2_{L^2}+\|e^{k+1}_\zeta\|^2_{L^2}\leq C\Delta{t}^4.\label{lemma 4.2 24}
\eena

Using the $H^2$ regularity to equation (\ref{Scheme3.1Numer Solu1}) and the boundedness of $\|\phi^{n}\|_{L^2}, \|q^{n}\|_{L^2}, \|\zeta^{n}\|_{L^2}$, we have
\begin{equation}
\begin{aligned}
\|\phi^{k+\frac{1}{2}}\|_{H^2}&\leq \|\phi^{k+\frac{1}{2}}\|_{L^2}+\|\frac{\phi^{k+1}-\phi^k}{\Delta{t}}\|_{L^2}+\|q^{k+\frac{1}{2}} \bar{g}^{k+\frac{1}{2}}\|_{L^2}+\|\zeta^{k+1/2}\|_{L^2}\\
&\leq \|\frac{e^{k+1}_\phi-e^k_\phi}{\Delta t}\|_{L^2}+\|\frac{\phi(t_{k+1})-\phi(t_{k})}{\Delta{t}}\|_{L^2}+\|\phi^{k+1}\|_{L^2}+\|\phi^{k}\|_{L^2}\\
&+\|q^{k+1}\|_{L^2}+\|q^{k}\|_{L^2}+\|\zeta^{k+1}\|_{L^2}+\|\zeta^{k}\|_{L^2}\leq C.\label{lemma 4.2 25}
\end{aligned}
\end{equation}
So
\bena
\|e^{k+\frac{1}{2}}_\phi\|_{H^2}\leq \|\phi(t_{k+1/2})\|_{H^2}+\|\phi^{k+1/2}\|_{H^2}\leq C.\label{lemma 4.2 26}$$
\eena
Finally, we have
\begin{equation}
\begin{aligned}
\|\phi^{k+\frac{1}{2}}\|_{L^\infty}\leq \|e^{k+\frac{1}{2}}_\phi\|_{L^\infty}+\|\phi(t_{k+1/2})\|_{L^\infty}
\leq \|e^{k+\frac{1}{2}}_\phi\|^{\frac{1}{2}}_{H^1}\|e^{k+\frac{1}{2}}_\phi\|^{\frac{1}{2}}_{H^2}+\|\phi(t_{k+1/2})\|_{L^\infty}
\leq C,\label{lemma 4.2 27}
\end{aligned}
\end{equation}
which implies
\bena
\|\phi^{k+1}\|_{L^\infty}= \|2\phi^{k+\frac{1}{2}}-\phi^{k}\|_{L^\infty}
\leq \|\phi^{k+\frac{1}{2}}\|_{L^\infty}+\|\phi^{k}\|_{L^\infty}
\leq C.\label{lemma 4.2 28}
\eena

By mathematical induction, $$\|\phi^n\|_{L^\infty}\leq C,n=0,1,\ldots,N.$$

\begin{thm}
For Scheme 3.1,  we have the following error estimates for $0\leq n\leq N$:
$$
\|e^{n}_\phi\|^2_{H^1}+\|e^{n}_q\|^2_{L^2}+\|e^{n}_\zeta\|^2_{L^2}\leq C\Delta{t}^4.
$$
\end{thm}
\noindent {\bf Proof.} Because $\|\phi^n\|_{L^\infty}(n=0,1,\ldots,N$) are bounded, the proof process is similar to (\ref{lemma4.2 1})-(\ref{lemma 4.2 24}) in Lemma 4.2 and is thus omitted.

\subsection{Error estimates for the Allen-Cahn  model of a Lagrange multiplier using EQ}

\noindent \indent Firstly, we note that
$$
  L=\frac{\int_{\Omega}[M\mu]d\bx}{\int_{\Omega}Md\bx}=\frac{1}{|\Omega|}\int_{\Omega}(-\gamma_{1}\Delta\phi+2\gamma_{2}\phi+2qg)\mathrm{d\bx}=\frac{1}{|\Omega|}\int_{\Omega}(2\gamma_{2}\phi+2qg)\mathrm{d\bx}.
$$
Then, Scheme 3.2 can be written as follows
\besl
&\frac{\phi^{n+1}-\phi^{n}}{\Delta t}=\gamma_{1}\Delta\phi^{n+1/2}-2\gamma_{2}\phi^{n+1/2}-2q^{n+1/2}\overline{g}^{n+1/2}+\frac{1}{\int_{\Omega}\mathrm{d\bx}}\int_{\Omega}(2\gamma_{2}\phi^{n+1/2}+2q^{n+1/2}\overline{g}^{n+1/2})\mathrm{d\bx},\label{Scheme3.2Numer Solu1}\\
&\frac{q^{n+1}-q^{n}}{\Delta t}=\overline{g}^{n+1/2}\frac{\phi^{n+1}-\phi^{n}}{\Delta t}.\label{Scheme3.2Numer Solu2}
\eesl

The error functions for $n\geq0$ satisfy the following difference equations
\bena\label{8}
\frac{e_{\phi}^{n+1}-e_{\phi}^{n}}{\Delta t}-\gamma_{1}\Delta e_{\phi}^{n+1/2}+2\gamma_{2}e_{\phi}^{n+1/2}+2q^{n+1/2}\overline{e_{g}}^{n+1/2}+2e_{q}^{n+1/2}\overline{g}^{n+1/2}\\
=\frac{1}{\int_{\Omega}\mathrm{d\bx}}\int_{\Omega}2\gamma_{2}e_{\phi}^{n+1/2}+2q^{n+1/2}\overline{e_{g}}^{n+1/2}+2e_{q}^{n+1/2}\overline{g}^{n+1/2}\mathrm{d\bx}+R_{\phi}^{n+1/2},
\eena

\bena\label{9}
\frac{e_{q}^{n+1}-e_{q}^{n}}{\Delta t}=\overline{g}^{n+1/2}\frac{e_{\phi}^{n+1}-e_{\phi}^{n}}{\Delta t}+\overline{e_{g}}^{n+1/2}\frac{\phi(t_{n+1})-\phi(t_{n})}{\Delta t}+R_{q}^{n+1/2}.
\eena

\begin{lem}
There exists a constant $C>0$ such that solution $\phi^n$ of Scheme 3.2 is bounded:
$$
\|\phi^n\|_{L^\infty}\leq C,n=0,1,\ldots,N.
$$
\end{lem}
\noindent {\bf Proof.} Taking inner product of (\ref{8}) with $e_\phi^{n+1}-e_\phi^n$, $2\Delta{t}e^{n+1/2}_\phi$, respectively, (\ref{9}) with $2e_q^{n+1}$, then adding them up and using inequalities, we obtain
\begin{equation}
\begin{aligned}
&\|e^{n+1}_\phi\|^2_{L^2}-\|e^n_\phi\|^2_{L^2}+\gamma_2\|{e}_\phi^{n+1}\|^2_{L^2}-\gamma_2\|{e}_\phi^{n}\|^2_{L^2}
+\frac{\gamma_1}{2}\|\nabla{e}^{n+1}_\phi\|^2_{L^2}-\frac{\gamma_1}{2}\|\nabla{e}^n_\phi\|^2_{L^2}+\|e^{n+1}_q\|^2_{L^2}-\|e^n_q\|^2_{L^2}\\
&\leq C\Delta{t}(\|e^{n+1}_\phi\|^2_{L^2}+\|e^n_\phi\|^2_{L^2}+\|e^{n-1}_\phi\|^2_{L^2}+\|\nabla{e}^{n+1}_\phi\|^2_{L^2}+\|\nabla{e}^{n}_\phi\|^2_{L^2}+\|\nabla{e}^{n-1}_\phi\|^2_{L^2}
+\|e^{n+1}_q\|^2_{L^2}+\|e^n_q\|^2_{L^2})+\Delta{t}O(\Delta{t}^4).
\end{aligned}
\end{equation}
The rest of the proof is similar to (\ref{lemma 4.2 20})-(\ref{lemma 4.2 28}) in the proof of Lemma 4.2 and is thus omitted.

\begin{thm}
For Scheme 3.2,  we have the following error estimates:
$$
\|e^{n}_\phi\|^2_{H^1}+\|e^{n}_q\|^2_{L^2}\leq C\Delta{t}^4,  0\leq n\leq N.
$$
\end{thm}
\noindent {\bf Proof.} Because $\|\phi^n\|_{L^\infty}(n=0,1,\ldots,N$) are bounded, the proof is similar to that in Lemma 4.3. So, we omit the details.

\subsection{Error estimates for the Allen-Cahn model of a penalizing potential using SAV}
\noindent \indent Scheme 3.3 can be written as
\besl
&\frac{\phi^{n+1}-\phi^n}{\Delta{t}}-\hs{\gamma_1 \Delta \phi}+2 \gamma_2 \hs{\phi} +2\hs{r} \ohs{s}+\sqrt{\eta}\zeta^{n+1/2}=0,\label{scheme 3.3 1}\\
&\frac{r^{n+1}-r^n}{\Delta t}= \int_\Omega\ohs{s}\frac{(\phi^{n+1}-\phi^n)}{\Delta t}\mathrm{d\bx} ,\label{scheme 3.3 2}\\
&\frac{\zeta^{n+1}-\zeta^n}{\Delta t}=\sqrt{\eta}\int_\Omega\frac{\phi^{n+1}-\phi^n}{\Delta t}\mathrm{d\bx}.\label{scheme 3.3 3}
\eesl
We deduce the error equations for $n\geq0$ as follows
\besl
&\frac{e_\phi^{n+1}-e_\phi^{n}}{\Delta{t}}-\gamma_1\Delta{e}_\phi^{n+1/2}+2\gamma_2e_\phi^{n+1/2}
+2r^{n+1/2}{\overline{e}_{s}^{n+1/2}}+2e_r^{n+1/2}\overline{s}^{n+1/2}+\sqrt{\eta}e^{n+1/2}_\zeta=R^{n+1/2}_\phi,\label{scheme 3.3 4}\\
&\frac{e_r^{n+1}-e_r^{n}}{\Delta{t}}=\int_\Omega{\overline{s}}^{n+1/2}\frac{e_\phi^{n+1}-e_\phi^{n}}{\Delta{t}}
+\overline{e}_{s}^{n+1/2}\frac{\phi(t_{n+1})-\phi(t_{n})}{\Delta{t}}\mathrm{d\bx}+R_r^{n+1/2},\label{scheme 3.3 5}\\
&\frac{e_\zeta^{n+1}-e_\zeta^{n}}{\Delta{t}}=\sqrt{\eta}\int_\Omega\frac{e_\phi^{n+1}-e_\phi^{n}}{\Delta{t}}\mathrm{d\bx}+R_\zeta^{n+1/2}.\label{scheme 3.3 6}
\eesl

\begin{lem}
There exists a constant $C>0$ such that the solution $\phi^n$ of Scheme 3.3 is bounded,
$$
\|\phi^n\|_{L^\infty}\leq C,n=0,1,\ldots,N.\label{scheme 3.3 7}
$$
\end{lem}
\noindent {\bf Proof.} Firstly, we note that $\|\phi^0\|_{L^\infty}\leq C$ is true  by definition. We assume $\|\phi^n\|_{L^\infty} (n=1, 2, \ldots, k)(k< N)$ are bounded. Then, we use mathematical induction to prove that $\|\phi^{k+1}\|_{L^\infty}$ is bounded.

We note that  $q(\phi)$ and $s(\phi)$  are continuous as functions of  $\phi$ and $\phi^n(n=1, 2, \ldots, k)$ and are thus bounded.
Taking inner product of (\ref{scheme 3.3 4}) with $e_\phi^{n+1}-e_\phi^n$, $2\Delta{t}e^{n+1/2}_\phi$, respectively, (\ref{scheme 3.3 5}) with $2e_r^{n+1}$ and (\ref{scheme 3.3 6}) with $2e_\zeta^{n+1}$, then adding them up, we obtain
\bena
\|{e}_\phi^{n+1}\|^2_{L^2}-\|{e}_\phi^{n}\|^2_{L^2}+\gamma_2(\|{e}_\phi^{n+1}\|^2_{L^2}-\|{e}_\phi^{n}\|^2_{L^2})+\frac{\gamma_1}{2}(\|\nabla{e}_\phi^{n+1}\|^2_{L^2}-\|\nabla{e}_\phi^{n}\|^2_{L^2})\\
+\|e^{n+1}_r\|^2_{L^2}-\|e^{n}_r\|^2_{L^2}+\|e^{n+1}_\zeta\|^2_{L^2}-\|e^{n}_\zeta\|^2_{L^2}+\frac{1}{\Delta{t}}\|e_\phi^{n+1}-e_\phi^{n}\|^2_{L^2}+\|e^{n+1}_r-e^{n}_r\|^2_{L^2}\\
+\|e^{n+1}_\zeta-e^{n}_\zeta\|^2_{L^2}+2\Delta{t}\gamma_1\|\nabla{e}^{n+1/2}_\phi\|^2_{L^2}+4\Delta{t}\gamma_2\|{e}^{n+1/2}_\phi\|^2_{L^2}\\
=(\int_\Omega{\overline{s(\phi)}}^{n+1/2}(e_\phi^{n+1}-e_\phi^{n})\mathrm{d\bx},2e^{n+1}_r)+(\int_\Omega\overline{e}_{s}^{n+1/2}(\phi(t_{n+1})-\phi(t_{n}))\mathrm{d\bx},2e^{n+1}_r)\\
+(\sqrt{\eta}\int_\Omega(e_\phi^{n+1}-e_\phi^{n})\mathrm{d\bx},2e^{n+1}_\zeta)-(2r^{n+1/2}{\overline{e}_{s}^{n+1/2}},e_\phi^{n+1}-e_\phi^n)\\
-(2e_r^{n+1/2}\overline{s(\phi)}^{n+1/2},e_\phi^{n+1}-e_\phi^n)-(\sqrt{\eta}e^{n+1/2}_\zeta,e_\phi^{n+1}-e_\phi^n)\\
-(2r^{n+1/2}{\overline{e}_{s}^{n+1/2}},2\Delta{t}e^{n+1/2}_\phi)-(2e_r^{n+1/2}\overline{{s}(\phi)}^{n+1/2},2\Delta{t}e^{n+1/2}_\phi)\\
-(\sqrt{\eta}e^{n+1/2}_\zeta,2\Delta{t}e^{n+1/2}_\phi)+(R^{n+1/2}_\phi,e_\phi^{n+1}-e_\phi^n)+\Delta{t}(R^{n+1/2}_r,2e^{n+1}_r)\\
+(R^{n+1/2}_\phi,2\Delta{t}e^{n+1/2}_\phi)+\Delta{t}(R^{n+1/2}_\zeta,2e^{n+1}_\zeta).\label{scheme 3.3 8}
\eena

Using Cauchy's inequality with $\epsilon$, we obtain
\bena
\label{scheme 3.3 9}
|(\int_\Omega{\overline{s}}^{n+1/2}(e_\phi^{n+1}-e_\phi^{n}){\mathrm{d\bx}},2e^{n+1}_r)|\\
\leq|(\|3{g}^n-{g}^{n-1}\|_{L^2}\|e_\phi^{n+1}-e_\phi^{n}\|_{L^2},e^{n+1}_r)|
\leq \frac{C}{4\Delta{t}\epsilon}\|e_\phi^{n+1}-e_\phi^{n}\|^2_{L^2}+C\Delta{t}\epsilon\|e^{n+1}_r\|^2_{L^2}.
\eena

By regularity assumption (\ref{regularity estimation}), Sobolev inequalities (\ref{Sobolev ineq 4.1}) and inequality (\ref{Lemma 4.1 inequality}), we get the following inequalities
\bena
\label{scheme 3.3 10}
|(\int_\Omega\overline{e}_{s}^{n+1/2}(\phi(t_{n+1})-\phi(t_{n})){ \mathrm{d\bx}},2e^{n+1}_r)|\\
\leq|(\|3e_{s}^n-e_{s}^{n-1}\|_{L^2}\|\phi(t_{n+1})-\phi(t_{n})\|_{L^2},e^{n+1}_r)|\\
\leq \|3e_{s}^n-e_{s}^{n-1}\|_{L^2}\|\phi(t_{n+1})-\phi(t_{n})\|_{L^2}\|e^{n+1}_r\|_{L^2}\\
\leq C\Delta{t}(\|3e_{s}^n\|^2_{L^4}+\|e_{s}^{n-1}\|^2_{L^4})+C\Delta{t}\|e^{n+1}_r\|^2_{L^2}\\
\leq C\Delta{t}(\|e^n_\phi\|^2_{L^2}+\|\nabla{e}^n_\phi\|^2_{L^2}+\|e^{n-1}_\phi\|^2_{L^2}+\|\nabla{e}^{n-1}_\phi\|^2_{L^2})+C\Delta{t}\|e^{n+1}_r\|^2_{L^2},
\eena
and
\bena
\label{scheme 3.3 11}
|(\sqrt{\eta}\int_\Omega{(e_\phi^{n+1}-e_\phi^{n})}\mathrm{d\bx},2e^{n+1}_\zeta)|\\
\leq\sqrt{\eta}(|\Omega|\|e_\phi^{n+1}-e_\phi^{n}\|_{L^2},2e^{n+1}_\zeta)|
\leq\frac{C}{4\Delta{t}\epsilon}\|e_\phi^{n+1}-e_\phi^n\|^2_{L^2}+C\Delta{t}\epsilon\|e_\zeta^{n+1}\|^2_{L^2}.
\eena

Analogous to the proof in Lemma 4.2, we obtain estimates for the  other items.
Choosing appropriate $\epsilon$ and combining the above inequalities (\ref{scheme 3.3 9})-(\ref{scheme 3.3 11}) and other items with (\ref{scheme 3.3 8}), we have
\bena
\|e^{n+1}_\phi\|^2_{L^2}+\gamma_2\|e^{n+1}_\phi\|^2_{L^2}+\frac{\gamma_1}{2}\|\nabla{e}^{n+1}_\phi\|^2_{L^2}+\|e^{n+1}_r\|^2_{L^2}+\|e^{n+1}_\zeta\|^2_{L^2}\\
-(\|e^{n}_\phi\|^2_{L^2}+\gamma_2\|e^{n}_\phi\|^2_{L^2}+\frac{\gamma_1}{2}\|\nabla{e}^{n}_\phi\|^2_{L^2}+\|e^{n}_r\|^2_{L^2}+\|e^{n}_\zeta\|^2_{L^2})\\
\leq C\Delta{t}(\|e^{n+1}_\phi\|^2_{L^2}+\|e^n_\phi\|^2_{L^2}+\|e^{n-1}_\phi\|^2_{L^2}+\|\nabla{e}^{n+1}_\phi\|^2_{L^2}+\|\nabla{e}^{n}_\phi\|^2_{L^2}\\
+\|\nabla{e}^{n-1}_\phi\|^2_{L^2}
+\|e^{n+1}_r\|^2_{L^2}+\|e^n_r\|^2_{L^2}+\|e^{n+1}_\zeta\|^2_{L^2}+\|e^n_\zeta\|^2_{L^2})\\
+\Delta{t}(\|R^{n+1/2}_\phi\|^2_{L^2}+\|R^{n+1/2}_r\|^2_{L^2}+\|R^{n+1/2}_\zeta\|^2_{L^2}).
\eena

The rest of the  proof  is the same as (\ref{lemma 4.2 20})-(\ref{lemma 4.2 28}) in Lemma 4.2. So, we skip the details and claim the following result
 $$\|\phi^n\|_{L^\infty}\leq C,n=0,1,\ldots,N.$$

\begin{thm}
For Scheme 3.3, we have the following error estimate:
$$
\|e^{n}_\phi\|^2_{H^1}+\|e^{n}_r\|^2_{L^2}+\|e^{n}_\zeta\|^2_{L^2}\leq C\Delta{t}^4,  0\leq n\leq N.
$$
\end{thm}
\noindent {\bf Proof.} Because $\|\phi^n\|_{L^\infty}(n=0,1,\ldots,N$) are bounded, the proof is similar to that in Lemma 4.4.

\subsection{Error estimates for the Allen-Cahn  model of a Lagrange multiplier using SAV}

\noindent \indent We rewrite Scheme 3.4 in the  following form,
\besl
&\frac{\phi^{n+1}-\phi^{n}}{\Delta t}=-[-\gamma_{1}\Delta\phi^{n+1/2}+2\gamma_{2}\phi^{n+1/2}+2r^{n+1/2}\overline{s}^{n+1/2}]+\frac{1}{\int_{\Omega}\mathrm{d\bx}}\int_{\Omega}(2\gamma_{2}\phi^{n+1/2}+2r^{n+1/2}\overline{s}^{n+1/2})\mathrm{d\bx},\label{scheme 3.4 1}\\
&\frac{r^{n+1}-r^{n}}{\Delta t}=\int_\Omega\overline{s}^{n+1/2}\frac{\phi^{n+1}-\phi^{n}}{\Delta t}\mathrm{d\bx}.\label{scheme 3.4 2}
\eesl
Then, we derive the equations for the error functions for $n\geq0$ as follows
\bena\label{scheme 3.3 12}
\frac{e_{\phi}^{n+1}-e_{\phi}^{n}}{\Delta t}-\gamma_{1}\Delta e_{\phi}^{n+1/2}+2\gamma_{2}e_{\phi}^{n+1/2}+2r^{n+1/2}\overline{e_{s}}^{n+1/2}+2e_{r}^{n+1/2}\overline{s}^{n+1/2}\\
=\frac{1}{\int_{\Omega}\mathrm{d\bx}}\int_{\Omega}2\gamma_{2}e_{\phi}^{n+1/2}+2r^{n+1/2}\overline{e_{s}}^{n+1/2}+2e_{r}^{n+1/2}\overline{s}^{n+1/2}\mathrm{d\bx}+R_{\phi}^{n+1/2},
\eena

\bena\label{scheme 3.3 13}
\frac{e_{r}^{n+1}-e_{r}^{n}}{\Delta t}=\int_\Omega\overline{s}^{n+1/2}\frac{e_{\phi}^{n+1}-e_{\phi}^{n}}{\Delta t}+\overline{e_{s}}^{n+1/2}\frac{\phi(t_{n+1})-\phi(t_{n})}{\Delta t}\mathrm{d\bx}+R_{r}^{n+1/2}.
\eena

\begin{lem}
There exists a constant $C>0$ such that the solution $\phi^n$ of Scheme 3.4 is bounded
$$
\|\phi^n\|_{L^\infty}\leq C,n=0,1,\ldots,N.
$$
\end{lem}
\noindent {\bf Proof.} Taking inner product of (\ref{scheme 3.3 12}) with $e_\phi^{n+1}-e_\phi^n$, $2\Delta{t}e^{n+1/2}_\phi$, respectively, (\ref{scheme 3.3 13}) with $2e_r^{n+1}$, then adding them up and using a series of inequalities, we obtain
\begin{equation}
\begin{aligned}
&\|e^{n+1}_\phi\|^2_{L^2}-\|e^n_\phi\|^2_{L^2}+\gamma_2\|{e}_\phi^{n+1}\|^2_{L^2}-\gamma_2\|{e}_\phi^{n}\|^2_{L^2}
+\frac{\gamma_1}{2}\|\nabla{e}^{n+1}_\phi\|^2_{L^2}-\frac{\gamma_1}{2}\|\nabla{e}^n_\phi\|^2_{L^2}+\|e^{n+1}_r\|^2_{L^2}-\|e^n_r\|^2_{L^2}\\
&\leq C\Delta{t}(\|e^{n+1}_\phi\|^2_{L^2}+\|e^n_\phi\|^2_{L^2}+\|e^{n-1}_\phi\|^2_{L^2}+\|\nabla{e}^{n+1}_\phi\|^2_{L^2}+\|\nabla{e}^{n}_\phi\|^2_{L^2}+\|\nabla{e}^{n-1}_\phi\|^2_{L^2}
+\|e^{n+1}_r\|^2_{L^2}+\|e^n_r\|^2_{L^2})+\Delta{t}O(\Delta{t}^4).
\end{aligned}
\end{equation}
The rest of the proof is similar to (\ref{lemma 4.2 20})-(\ref{lemma 4.2 28}) in that of  Lemma 4.2 and is thus omitted.

\begin{thm}
For Scheme 3.4, we have the following error estimates
$$
\|e^{n}_\phi\|^2_{H^1}+\|e^{n}_r\|^2_{L^2}\leq C\Delta{t}^4, 0\leq n\leq N.
$$
\end{thm}
\noindent {\bf Proof.} Because $\|\phi^n\|_{L^\infty}(n=0,1,\ldots,N$) are bounded, the proof is similar to that in Lemma 4.5 and thus omitted.

\subsection{{Error estimates of fully discrete numerical schemes}}

\noindent \indent In this subsection, we omit the space symbol $i, j$ for simplicity. Then, fully discrete Scheme 3.5 can be written as follows
\besl
&\frac{\phi^{n+1}-\phi^{n}}{\Delta t}=\gamma_{1}\Delta_h\phi^{n+1/2}-2\gamma_{2}\phi^{n+1/2}-2q^{n+1/2}\overline{g}^{n+1/2}+\frac{1}{\int_{\Omega}\mathrm{d\bx}}[1\star (2\gamma_{2}\phi^{n+1/2}+2q^{n+1/2}\overline{g}^{n+1/2})],\label{Scheme3.5Numer Solu1}\\
&\frac{q^{n+1}-q^{n}}{\Delta t}=\overline{g}^{n+1/2}\frac{\phi^{n+1}-\phi^{n}}{\Delta t}.\label{Scheme3.5Numer Solu2}
\eesl
The error functions for $n\geq0$ satisfy the following difference equations
\bena\label{Scheme3.5error function1}
\frac{e_{\phi}^{n+1}-e_{\phi}^{n}}{\Delta t}-\gamma_{1}\Delta_h e_{\phi}^{n+1/2}+2\gamma_{2}e_{\phi}^{n+1/2}+2q^{n+1/2}\overline{e_{g}}^{n+1/2}+2e_{q}^{n+1/2}\overline{g}^{n+1/2}\\
=\frac{1}{\int_{\Omega}\mathrm{d\bx}}[1\star (2\gamma_{2}e_{\phi}^{n+1/2}+2q^{n+1/2}\overline{e_{g}}^{n+1/2}+2e_{q}^{n+1/2}\overline{g}^{n+1/2})]
+\tau_{\phi}^{n+1/2},
\eena
\bena\label{Scheme3.5error function2}
\frac{e_{q}^{n+1}-e_{q}^{n}}{\Delta t}=\overline{g}^{n+1/2}\frac{e_{\phi}^{n+1}-e_{\phi}^{n}}{\Delta t}+\overline{e_{g}}^{n+1/2}\frac{\phi(t_{n+1})-\phi(t_{n})}{\Delta t}+\tau_{q}^{n+1/2},
\eena
where the corresponding truncation error
$$|\tau^{n+1/2}_\phi|\leq C_1({h}^2+{\Delta t}^2), |\tau^{n+1/2}_q|\leq C_2{\Delta t}^2.$$

\begin{lem}
There exists a constant $C>0$ such that solution $\phi^n$ of Scheme 3.5 is bounded:
$$\|\phi^n\|_\infty\leq C,n=0,1,\ldots,N.$$
\end{lem}
\noindent {\bf Proof.} Firstly, we note that  $\|\phi^0\|_\infty\leq C$ is true by definition. We assume $\|\phi^n\|_\infty (n=1, 2, \ldots, k)(k< N)$ are bounded. Next, we use mathematical induction to prove $\|\phi^{k+1}\|_\infty$ is bounded.

We note that $q(\phi)$ and $g(\phi)$ are continuous as functions of  $\phi$ and $\phi^n(n=1, 2, \ldots, k)$ and thus are bounded.
Taking discrete inner product of (\ref{Scheme3.5error function1}) with $e_\phi^{n+1}-e_\phi^n$ and $2\Delta{t}e^{n+1/2}_\phi$, respectively, summing over $i$ and $j$,  we obtain
\bena
\frac{1}{\Delta{t}}\|e_\phi^{n+1}-e_\phi^{n}\|^2_{2}+\frac{\gamma_1}{2}(\|\nabla_h{e}_\phi^{n+1}\|^2_{2}-\|\nabla_h{e}_\phi^{n}\|^2_{2})
+\gamma_2(\|{e}_\phi^{n+1}\|^2_{2}-\|{e}_\phi^{n}\|^2_{2})+(2q^{n+1/2}{\overline{e}_{g}^{n+1/2}},e_\phi^{n+1}-e_\phi^n)\\
+(2e_q^{n+1/2}\overline{g(\phi)}^{n+1/2},e_\phi^{n+1}-e_\phi^n)-([1\star (2\gamma_{2}e_{\phi}^{n+1/2}+2q^{n+1/2}\overline{e_{g}}^{n+1/2}+2e_{q}^{n+1/2}\overline{g}^{n+1/2})],e_\phi^{n+1}-e_\phi^n)\\
=(\tau^{n+1/2}_\phi,e_\phi^{n+1}-e_\phi^n)\label{lemma6 1}
\eena
and
\bena
\|{e}_\phi^{n+1}\|^2_{2}-\|{e}_\phi^{n}\|^2_{2}+2\Delta{t}\gamma_1\|\nabla_h{e}^{n+1/2}_\phi\|^2_{2}+4\Delta{t}\gamma_2\|{e}^{n+1/2}_\phi\|^2_{2}+(2q^{n+1/2}{\overline{e}_{g}^{n+1/2}},2\Delta{t}e^{n+1/2}_\phi)\\
+(2e_q^{n+1/2}\overline{g(\phi)}^{n+1/2},2\Delta{t}e^{n+1/2}_\phi)-([1\star (2\gamma_{2}e_{\phi}^{n+1/2}+2q^{n+1/2}\overline{e_{g}}^{n+1/2}+2e_{q}^{n+1/2}\overline{g}^{n+1/2})],2\Delta{t}e^{n+1/2}_\phi)\\
=(\tau^{n+1/2}_\phi,2\Delta{t}e^{n+1/2}_\phi).\label{lemma6 2}
\eena

Taking discrete inner product of (\ref{Scheme3.5error function2}) with $2e_q^{n+1}$, summing over $i$ and $j$, we have
\bena
\|e^{n+1}_q\|^2_{2}-\|e^{n}_q\|^2_{2}+\|e^{n+1}_q-e^{n}_q\|^2_{2}\\
=({\overline{g}}^{n+1/2}(e_\phi^{n+1}-e_\phi^{n}),2e^{n+1}_q)+(\overline{e}_{g}^{n+1/2}(\phi(t_{n+1})-\phi(t_{n})),2e^{n+1}_q)
+\Delta{t}(\tau^{n+1/2}_q,2e^{n+1}_q).\label{lemma6 3}
\eena
Combining (\ref{lemma6 1})-(\ref{lemma6 3}),we have
\bena
\|{e}_\phi^{n+1}\|^2_{2}-\|{e}_\phi^{n}\|^2_{2}+\gamma_2(\|{e}_\phi^{n+1}\|^2_{2}-\|{e}_\phi^{n}\|^2_{2})
+\frac{\gamma_1}{2}(\|\nabla_h{e}_\phi^{n+1}\|^2_{2}-\|\nabla_h{e}_\phi^{n}\|^2_{2})+\|e^{n+1}_q\|^2_{2}-\|e^{n}_q\|^2_{2}\\
+\frac{1}{\Delta{t}}\|e_\phi^{n+1}-e_\phi^{n}\|^2_{2}+\|e^{n+1}_q-e^{n}_q\|^2_{2}+
+2\Delta{t}\gamma_1\|\nabla_h{e}^{n+1/2}_\phi\|^2_{2}+4\Delta{t}\gamma_2\|{e}^{n+1/2}_\phi\|^2_{2}\\
=({\overline{g(\phi)}}^{n+1/2}(e_\phi^{n+1}-e_\phi^{n}),2e^{n+1}_q)+(\overline{e}_{g}^{n+1/2}(\phi(t_{n+1})-\phi(t_{n})),2e^{n+1}_q)\\
-(2q^{n+1/2}{\overline{e}_{g}^{n+1/2}},e_\phi^{n+1}-e_\phi^n)-(2e_q^{n+1/2}\overline{g(\phi)}^{n+1/2},e_\phi^{n+1}-e_\phi^n)\\
-(2q^{n+1/2}{\overline{e}_{g}^{n+1/2}},2\Delta{t}e^{n+1/2}_\phi)-(2e_q^{n+1/2}\overline{{g}(\phi)}^{n+1/2},2\Delta{t}e^{n+1/2}_\phi)\\
+([1\star (2\gamma_{2}e_{\phi}^{n+1/2}+2q^{n+1/2}\overline{e_{g}}^{n+1/2}+2e_{q}^{n+1/2}\overline{g}^{n+1/2})],e_\phi^{n+1}-e_\phi^n)\\
+([1\star (2\gamma_{2}e_{\phi}^{n+1/2}+2q^{n+1/2}\overline{e_{g}}^{n+1/2}+2e_{q}^{n+1/2}\overline{g}^{n+1/2})],2\Delta{t}e^{n+1/2}_\phi)\\
+(\tau^{n+1/2}_\phi,e_\phi^{n+1}-e_\phi^n)+\Delta{t}(\tau^{n+1/2}_q,2e^{n+1}_q)+(\tau^{n+1/2}_\phi,2\Delta{t}e^{n+1/2}_\phi).\label{lemma6 4}
\eena

Similar to the proof of Lemma 4.2, we apply the corresponding inequalities to have
\bena
\label{scheme6 5}
|([1\star (2\gamma_{2}e_{\phi}^{n+1/2}+2q^{n+1/2}\overline{e_{g}}^{n+1/2}+2e_{q}^{n+1/2}\overline{g}^{n+1/2})],e_\phi^{n+1}-e_\phi^n)|\\
\leq|(\|2\gamma_{2}e_{\phi}^{n+1/2}+2q^{n+1/2}\overline{e_{g}}^{n+1/2}+2e_{q}^{n+1/2}\overline{g}^{n+1/2}\|_{2},e_\phi^{n+1}-e_\phi^n)|\\
\leq C\Delta{t}\epsilon(\|2\gamma_{2}e_{\phi}^{n+1/2}+2q^{n+1/2}\overline{e_{g}}^{n+1/2}+2e_{q}^{n+1/2}\overline{g}^{n+1/2}\|_{2}^2)+\frac{C}{4\Delta{t}\epsilon}(\|e_\phi^{n+1}-e_\phi^n\|^2_2)\\
\leq C\Delta{t}\epsilon(\|e_\phi^{n+1}\|^2_{2}+\|e_\phi^{n}\|^2_{2}+\|e^n_\phi\|^2_{2}
+\|\nabla_h{e}^n_\phi\|^2_{2}+\|e^{n-1}_\phi\|^2_{2}+\|\nabla_h{e}^{n-1}_\phi\|^2_{2}\\
+\|e^{n+1}_q\|^2_{2}+\|{e}^n_q\|^2_{2})+\frac{C}{4\Delta{t}\epsilon}(\|e_\phi^{n+1}-e_\phi^n\|^2_2)
\eena
and
\bena
\label{scheme6 6}
|([1\star (2\gamma_{2}e_{\phi}^{n+1/2}+2q^{n+1/2}\overline{e_{g}}^{n+1/2}+2e_{q}^{n+1/2}\overline{g}^{n+1/2})],2\Delta{t}e^{n+1/2}_\phi)|\\
\leq \Delta{t}(\|2\gamma_{2}e_{\phi}^{n+1/2}+2q^{n+1/2}\overline{e_{g}}^{n+1/2}+2e_{q}^{n+1/2}\overline{g}^{n+1/2}\|^2_2+\|e^{n+1/2}_\phi\|^2_2)\\
\leq C\Delta{t}(\|e_\phi^{n+1}\|^2_{2}+\|e_\phi^{n}\|^2_{2}+\|e^n_\phi\|^2_{2}
+\|\nabla_h{e}^n_\phi\|^2_{2}+\|e^{n-1}_\phi\|^2_{2}+\|\nabla_h{e}^{n-1}_\phi\|^2_{2}\\
+\|e^{n+1}_q\|^2_{2}+\|{e}^n_q\|^2_{2}).
\eena

Thus, we have the following result
\bena
\|e^{n+1}_\phi\|^2_{2}+\gamma_2\|e^{n+1}_\phi\|^2_{2}+\frac{\gamma_1}{2}\|\nabla_h{e}^{n+1}_\phi\|^2_{2}+\|e^{n+1}_q\|^2_{2}-(\|e^{n}_\phi\|^2_{2}+\gamma_2\|e^{n}_\phi\|^2_{2}+\frac{\gamma_1}{2}\|\nabla_h{e}^{n}_\phi\|^2_{2}+\|e^{n}_q\|^2_{2})\\
\leq C\Delta{t}(\|e^{n+1}_\phi\|^2_{2}+\|e^n_\phi\|^2_{2}+\|e^{n-1}_\phi\|^2_{2}+\|\nabla_h{e}^{n}_\phi\|^2_{2}+\|\nabla_h{e}^{n-1}_\phi\|^2_{2}
+\|e^{n+1}_q\|^2_{2}+\|e^n_q\|^2_{2})\\
+\Delta{t}(\|\tau^{n+1/2}_\phi\|^2_{2}+\|\tau^{n+1/2}_q\|^2_{2})\\
\leq C\Delta{t}(\|e^{n+1}_\phi\|^2_{2}+\|e^n_\phi\|^2_{2}+\|e^{n-1}_\phi\|^2_{2}+\|\nabla_h{e}^{n+1}_\phi\|^2_{2}+\|\nabla_h{e}^{n}_\phi\|^2_{2}+\|\nabla_h{e}^{n-1}_\phi\|^2_{2}
+\|e^{n+1}_q\|^2_{2}+\|e^n_q\|^2_{2})\\
+\Delta{t}(\|\tau^{n+1/2}_\phi\|^2_{2}+\|\tau^{n+1/2}_q\|^2_{2}).\label{lemma6 7}
\eena

Summing (\ref{lemma6 7}) over for time steps 1 to m
and applying  the Gronwall's inequality and discrete Sobolev inequalities
\bena
\|f\|^2_{4}\leq C\|f\|_{2}\|f\|_{H_h^1},
\|f\|^2_{\infty}\leq C\|f\|_{H_h^1}\|f\|_{H_h^2},\label{Sobolev dis ineq 4.1}
\eena
 we arrive at
\bena
\|e^{k+1}_\phi\|^2_{H^1_h}+\|e^{k+1}_q\|^2_{2}\leq C({h}^2+{\Delta t}^2)^2.\label{lemma6 8}
\eena

The rest of the proof is then similar to (\ref{lemma 4.2 25})-(\ref{lemma 4.2 28}) in the proof of Lemma 4.2. So, we have
\ben
\|\phi^n\|_\infty\leq C,n=0,1,\ldots,N.
\een

\begin{thm}
For numerical solutions of Scheme 3.5,  we have the following fully discrete error estimate for $0\leq n\leq N$:
$$
\|e^{n}_\phi\|^2_{H^1_h}+\|e^{n}_q\|^2_{2}\leq C(h^2+\Delta{t}^2)^2.
$$
\end{thm}
\noindent {\bf Proof.} Because $\|\phi^n\|_\infty(n=0,1,\ldots,N$) are bounded, the proof  is similar to (\ref{lemma6 1})-(\ref{lemma6 8}) in the proof of Lemma 4.6 and is thus omitted.

\section{Numerical Results}

\noindent \indent In this section, we conduct some numerical experiments to show the accuracy and usefulness of the numerical schemes.
In the numerical experiments, we use square domain: $\Omega=[0,1]\times[0,1]$ and  model parameter values $\eta=10^4, C_0=10^4$ wherever relevant.

\subsection{Mesh refinement}

\noindent \indent We first use the initial condition
$$\phi(x,y,0)=\frac{1}{2}+\frac{1}{2}cos(4\pi x)cos(4\pi y)$$
for the mesh refinement test and set the  parameter values as $\gamma_1=0.2, \gamma_2=10$ and $M=1\times 10^{-3}$. Since we do not have the exact solution, we use the difference between results on successive coarse and finer grids to evaluate the numerical errors. We output the errors at time $T=1.0$. The mesh refinement test results are summarized in  Table 5.1-5.4 for temporal mesh refinement at fixed $h=1/256$ and in Table 5.5 for spatial mesh refinement at fixed $\Delta t=1.0\times10^{-4}$. Second order convergence rates are observed in all the tests.

\begin{table}[htbp]
\centering
\caption{Convergence rates in time of Scheme 3.1}
\begin{tabular}{ccccccc}
\hline
~& Coarse $\Delta{t}$~    ~& Fine $\Delta{t}$~       ~&$L^2$-Error of Scheme      ~&Rates           ~& $H^1$-Error of Scheme        ~&Rates          \\ \hline
~&$0.1$~                  ~&$0.05$~          & 7.28e-07            & -                    & 5.64e-04               &  -              \\
~&$0.05$~                  ~&$0.025$~          & 1.81e-07            & 2.003523               & 1.39e-04               &2.016886              \\
~&$0.025$~                  ~&$0.0025$~         & 4.53e-08            & 2.001735               & 3.46e-05               &2.008339                 \\
~&$0.0025$~                 ~&$0.00125$~         & 1.13e-08            & 2.000862               & 8.63e-06               &2.004144                   \\
~&$0.00125$~                ~&$0.000625$~          & 2.83e-09           & 2.00043               & 2.15e-06               &2.002065                     \\ \hline
\end{tabular}
\end{table}

\begin{table}[htbp]
\centering
\caption{Convergence rates in time of Scheme 3.2}
\begin{tabular}{ccccccc}
\hline
~& Coarse $\Delta{t}$~    ~& Fine $\Delta{t}$~          ~&$L^2$-Error of Scheme      ~&Rates           ~& $H^1$-Error of Scheme        ~&Rates      \\ \hline
~&$0.1$~                    ~&$0.05$~                            & 7.28e-07            & -                    & 5.64e-04               &  -              \\
~&$0.05$~                   ~&$0.025$~                            & 1.81e-07            & 2.003523               & 1.39e-04               &2.016886      \\
~&$0.025$~                  ~&$0.0025$~                           & 4.53e-08            & 2.001735               & 3.46e-05               &2.008339         \\
~&$0.0025$~                 ~&$0.00125$~                           & 1.13e-08            & 2.000861               & 8.63e-06               &2.004143       \\
~&$0.00125$~                ~&$0.000625$~                          & 2.83e-09           & 2.000424               & 2.15e-06               &2.002062   \\ \hline
\end{tabular}
\end{table}

\begin{table}[htbp]
\centering
\caption{Convergence rates in time of Scheme 3.3}
\begin{tabular}{ccccccc}
\hline
~& Coarse $\Delta{t}$~    ~& Fine $\Delta{t}$~          ~&$L^2$-Error of Scheme      ~&Rates           ~& $H^1$-Error of Scheme        ~&Rates     \\ \hline
~&$0.1$~     ~&$0.05$~                            & 7.27e-07            & -                    & 5.64e-04               &  -        \\
~&$0.05$~    ~&$0.025$~                            & 1.81e-07            & 2.003524               & 1.39e-04               &2.01689     \\
~&$0.025$~    ~&$0.0025$~                           & 4.53e-08            & 2.001735               & 3.46e-05               &2.008341      \\
~&$0.0025$~    ~&$0.00125$~                           & 1.13e-08            & 2.000861               & 8.63e-06               &2.004145      \\
~&$0.00125$~    ~&$0.000625$~                          & 2.83e-09           & 2.000434               & 2.15e-06               &2.002067     \\ \hline
\end{tabular}
\end{table}

\begin{table}[htbp]
\centering
\caption{Convergence rates in time of Scheme 3.4}
\begin{tabular}{ccccccc}
\hline
~& Coarse $\Delta{t}$~      ~& Fine $\Delta{t}$~          ~&$L^2$-Error of Scheme      ~&Rates           ~& $H^1$-Error of Scheme        ~&Rates  \\ \hline
~&$0.1$~         ~&$0.05$~                            & 7.27e-07            & -                    & 5.64e-04               &  -    \\
~&$0.05$~         ~&$0.025$~                            & 1.81e-07            & 2.003524               & 1.39e-04               &2.01689   \\
~&$0.025$~          ~&$0.0025$~                           & 4.53e-08            & 2.001735               & 3.46e-05               &2.008341   \\
~&$0.0025$~         ~&$0.00125$~                           & 1.13e-08            & 2.000861               & 8.63e-06               &2.004145   \\
~&$0.00125$~          ~&$0.000625$~                          & 2.83e-09           & 2.000433               & 2.15e-06               &2.002066    \\ \hline
\end{tabular}
\end{table}

\begin{table}[htbp]
\centering
\caption{Convergence rates in space of Scheme 3.5}
\begin{tabular}{ccccccc}
\hline
~& Coarse $h$~     ~& Fine $h$~          ~&$L^2$-Error of Scheme      ~&Rates           ~& $H^1$-Error of Scheme        ~&Rates \\ \hline
~&$0.125$~        ~&$0.0625$~                            & 0.05066            & -               & 9.964354               & -     \\
~&$0.0625$~       ~&$0.03125$~                           & 0.013207            & 1.93954               & 2.830995               &1.815467 \\
~&$0.03125$~      ~&$0.015625$~                           & 0.003403            & 1.956586               & 0.736601               &1.942353  \\
~&$0.015625$~      ~&$0.0078125$~                          & 9.03e-04           & 1.913851               & 0.195949               &1.910403    \\
~&$0.0078125$~      ~&$0.00390625$~                          & 2.07e-04           & 2.123257               & 0.044999               &2.122508    \\ \hline
\end{tabular}
\end{table}

\subsection{Merging drops}

\noindent \indent To test the properties of volume-conservation and energy dissipation for the Allen-Cahn equations with nonlocal constraints, we conduct a numerical experiment in which four drops merge into a large one. We set the initial condition as follows
\begin {equation*}
\left\{
\begin{array}{lr}
1,& r_1\leq0.2-\delta \quad \text{or}\quad r_2\leq0.2-\delta \quad \text{or}\quad r_3\leq0.2-\delta \quad \text{or}\quad r_4\leq0.2-\delta,\\
tanh(\frac{0.2+\delta-r_1}{\delta}),& 0.2-\delta<r_1<0.2+\delta,\\
tanh(\frac{0.2+\delta-r_2}{\delta}),& 0.2-\delta<r_2<0.2+\delta,\\
tanh(\frac{0.2+\delta-r_3}{\delta}),& 0.2-\delta<r_3<0.2+\delta,\\
tanh(\frac{0.2+\delta-r_4}{\delta}),& 0.2-\delta<r_4<0.2+\delta,\\
0,& \text {other},
\end{array}
\right.
\end{equation*}
where $\delta=0.01$, $r_1=\sqrt{(x-0.3+\delta)^2+(y-0.3+\delta)^2}$, $r_2=\sqrt{(x-0.7-\delta)^2+(y-0.3+\delta)^2}$,\\
 $r_3=\sqrt{(x-0.3+\delta)^2+(y-0.7-\delta)^2}$, $r_4=\sqrt{(x-0.7-\delta)^2+(y-0.7-\delta)^2}$.

In this example, we set the  parameter values as $\gamma_1=0.02, \gamma_2=100$, and $M=1$. We use  $\Delta t=1.0\times10^{-4}$ and $N_x=N_y=256$  in the simulation. The numerical results calculated by the EQ schemes and SAV schemes for the same model are consistent (Figure 5.1 and 5.2). It shows that the classical Allen-Cahn model does not  conserve the volume, which leads to disappearance of the drops after sometime (Figure5.1A and Figure5.2A). However, the model with a penalizing potential and the model with a lagrange multiplier  conserve the volume pretty well ( Figure5.1B and C,  Figure5.2B and C, Figure5.3B). Concerning energy dissipation, the energy curve of the Allen-Cahn model will decay to zero until the drop disappears (Figure5.3A), whereas the energy curves of both modified models with  nonlocal constraints decrease with time until a steady state plateau is reached, which means the four initial droplets have merged into one large droplet with its volume equal to the combined volume of the four small drops (Figure5.1B and C, Figure5.2B and C, Figure5.3A and B). Both volume conservation and energy dissipation in the schemes are clearly demonstrated by the example.

\begin{figure}[htbp]
  \centering
  \includegraphics[width=.7 \textwidth]{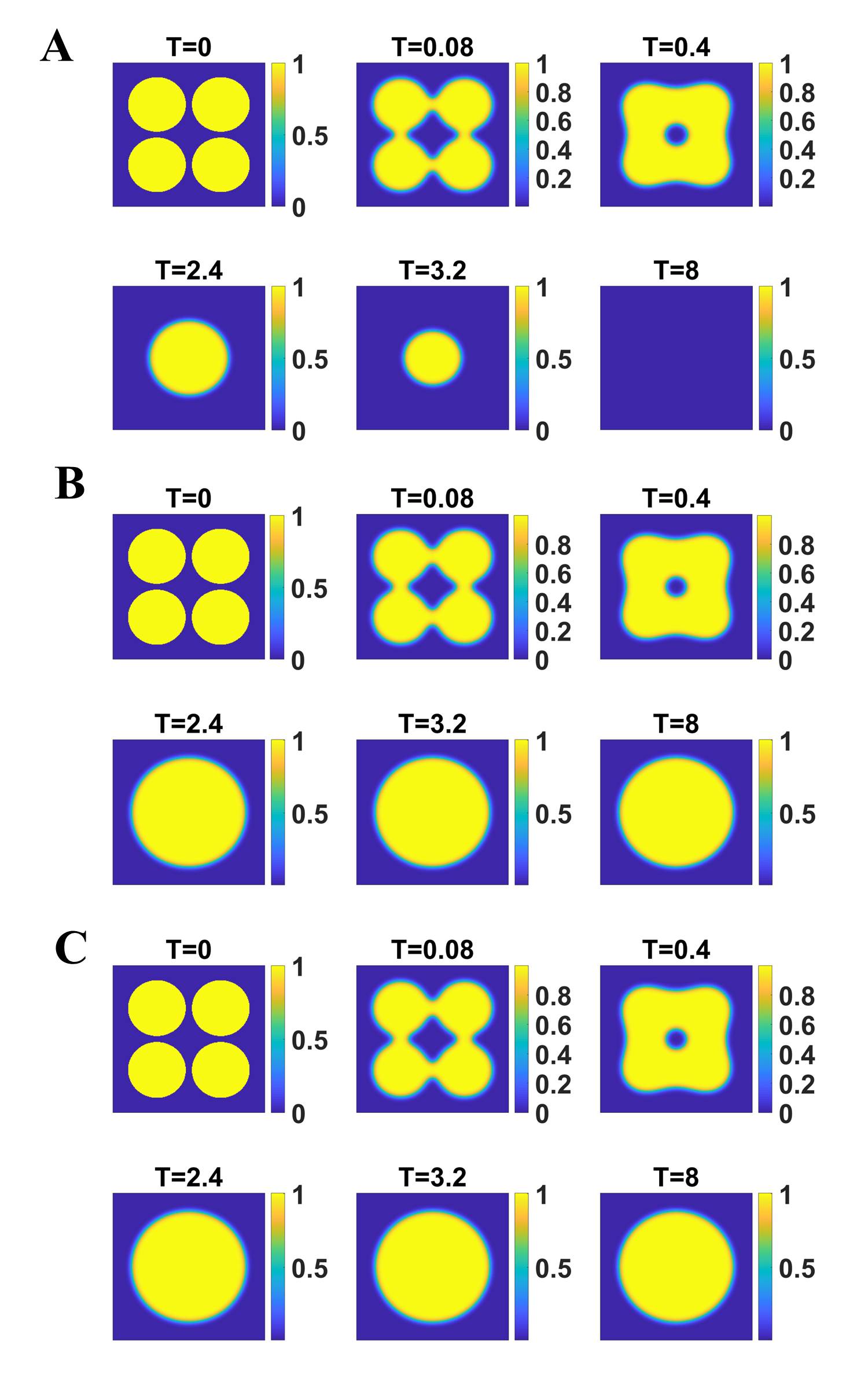}
  \caption{ The snapshots of numerical simulations of four drop merging experiments at $T=0, 0.08, 0.4, 2.4, 3.2, 8$ for (A) the Allen-Cahn model, (B) the Allen-Cahn model with a penalizing potential and (C) the Allen-Cahn model with a Lagrange multiplier using EQ methods, respectively. Clearly, the Allen-Cahn model gives the wrong result in the long time simulation while the modified Allen-Cahn models with nonlocal constraints preserve the volume of the drops and dissipate the total energy as well.}
\end{figure}
\clearpage
\begin{figure}[htbp]
  \centering
  \includegraphics[width=.7 \textwidth]{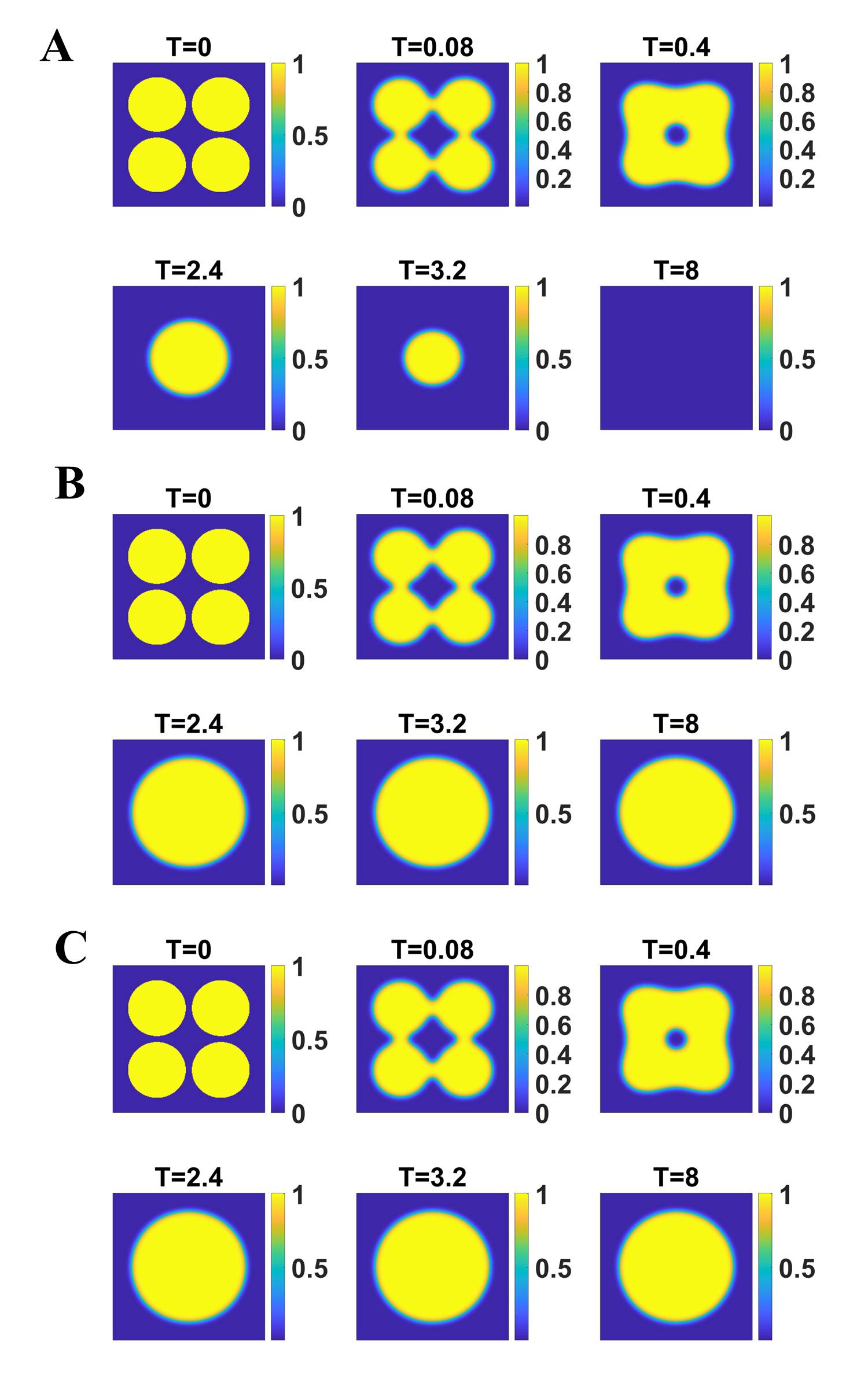}
  \caption{The snapshots of numerical simulations of four drop merging at $T=0, 0.08, 0.4, 2.4, 3.2, 8$ for (A) the Allen-Cahn model, (B) the Allen-Cahn model with a penalizing potential  and (C) the Allen-Cahn model with a Lagrange multiplier using SAV methods, respectively. Clearly, the Allen-Cahn model gives the wrong result in long time simulations while the modified Allen-Cahn models with nonlocal constraints preserve the volume of the drops and dissipate the total energy as well.}
\end{figure}
\clearpage
\begin{figure}[htbp]
  \centering
  \includegraphics[width=1 \textwidth]{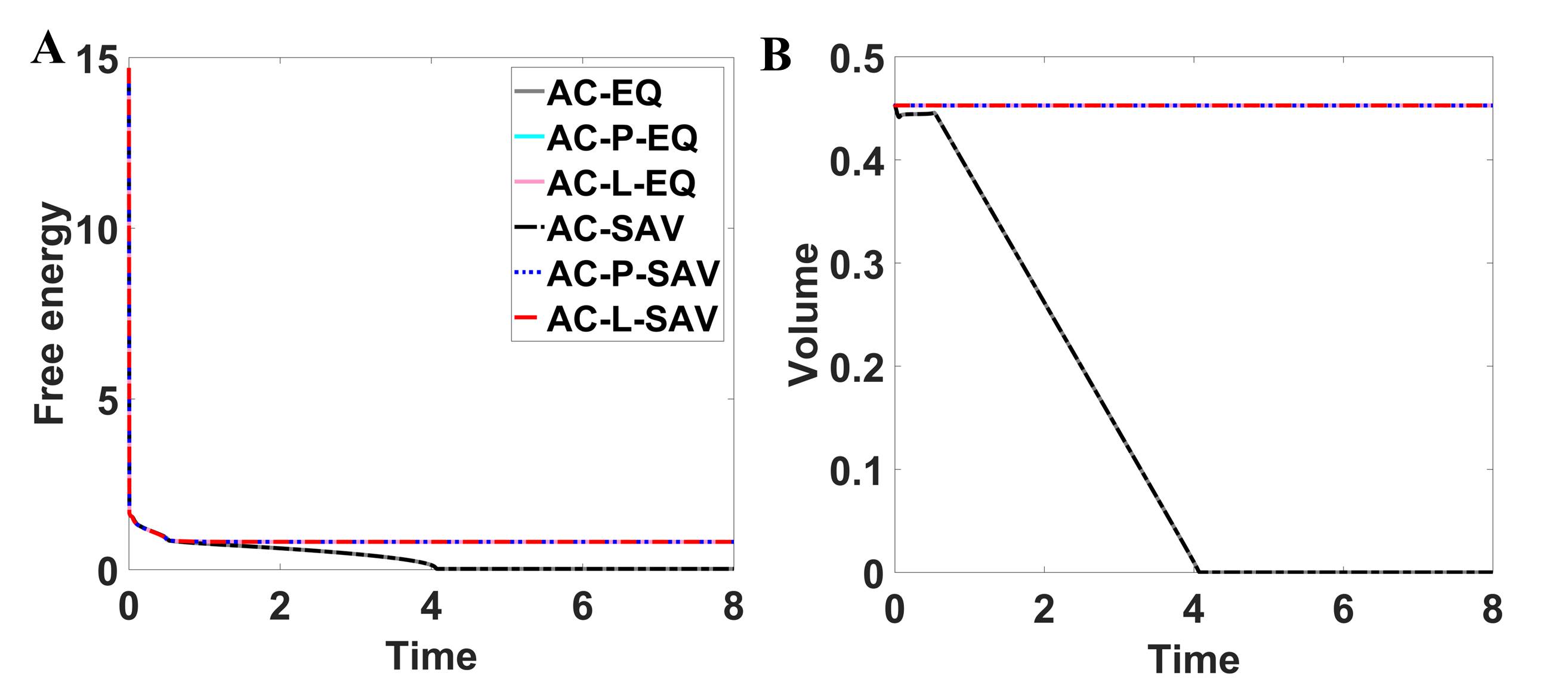}
  \caption{Time evolution of the volume and the free energy of the Allen-Cahn model and the Allen-Cahn models with nonlocal constraints. In short time, the Allen-Cahn models give similar results; whereas in long time, the Allen-Cahn model yields a zero volume and zero free energy while the modified Allen-Cahn models with nonlocal constraints preserve the volume and reaches a nearly steady state with a non-zero free energy. AC-\{EQ,SAV\}: the Allen-Cahn equation discretized using \{EQ, SAV\} methods. AC-P-\{EQ,SAV\}: the Allen-Cahn equation with a penalizing potential discretized using \{EQ,SAV\} methods. AC-L-\{EQ,SAV\}: the Allen-Cahn equation with a Lagrange multiplier discretized using \{EQ, SAV\} methods. }
\end{figure}

\section{Conclusion}

\noindent \indent We have developed two pairs of numerical schemes for  two    Allen-Cahn equations with nonlocal constraints using EQ and SAV methods, respectively. The schemes are shown as unconditionally  energy stable and volume-preserving.  The linear systems resulting from the schemes are uniquely solvable. We then carry out a series of error estimates for the four semi-discrete numerical schemes and one representative fully discrete scheme, and verify their second-order convergence rates in time {and space}, respectively. Numerical simulations are conducted to illustrate the accuracy of the proposed numerical schemes and their conservative properties. Under the volume-preserving constraints, both modified Allen-Cahn  models with nonlocal constraints preserve the volume pretty well in the numerical experiments. Therefore, the Allen-Cahn models with nonlocal constraints can be used as an effective alternative to the Cahn-Hilliard model when modeling multi-component material systems where phase volume and energy dissipation are important quantities to preserve.

\section*{Acknowledgements}
 %Jun Li's work is supported by the National Natural Science Foundation of China (Grant No. 11301287).
% Jia Zhao and Qi Wang are partially supported by NSF-DMS-1200487 and  DMS-1517347 award.
 Qi Wang's research is partially supported by NSF awards DMS-1517347,  DMS-1815921 and OIA-1655740, and NSFC awards \#11571032, \#91630207 and NSAF-U1530401.
%\begin{appendices}
%\section{}
%\subsection{}
%\section{}
%\subsection{}
%\end{appendices}
\bibliographystyle{plain}
\bibliography{reference}
\end{document}